\documentclass[11pt]{article}

\usepackage{amsmath,amsthm,amssymb,amscd}

\usepackage{graphicx}
\usepackage{epstopdf}
\usepackage{color}
\usepackage{url}
\usepackage{algorithm}
\usepackage{algorithmic}

\numberwithin{equation}{section}

\theoremstyle{definition}\newtheorem{definition}{Definition}
\theoremstyle{plain}\newtheorem{theorem}[definition]{Theorem}
\theoremstyle{definition}
\theoremstyle{definition}
\theoremstyle{definition}
\theoremstyle{definition}
\theoremstyle{definition}
\theoremstyle{definition}\newtheorem{lemma}[definition]{Lemma}
\theoremstyle{definition}

\newcommand{\cR}{{\mathcal R}}

\DeclareMathOperator*{\argmin}{argmin}
\newcommand{\R}{\mathbb{R}}
\newcommand{\N}{\mathbb{N}}

\newcommand{\xad}{x_\alpha^\delta}
\newcommand{\xd}{x^\dag}
\newcommand{\xa}{x_\alpha}

\newcommand{\asa}{A^\ast A}
\newcommand{\basa}{(A^\ast A)}
\newcommand{\bigo}{\mathcal{O}}

\begin{document}

\title{Estimating solution smoothness and data noise with Tikhonov regularization}

\author{\textsc{Daniel Gerth, Ronny Ramlau}
}


\maketitle

\begin{abstract}
A main drawback of classical Tikhonov regularization is that often the parameters required to apply theoretical results, e.g., the smoothness of the sought-after solution and the noise level, are unknown in practice. In this paper we investigate in new detail the residuals in Tikhonov regularization viewed as functions of the regularization parameter. We show that the residual carries, with some restrictions, the information on both the unknown solution and the noise level. By calculating approximate solutions for a large range of regularization parameters, we can extract both parameters from the residual given only one set of noisy data and the forward operator. The smoothness in the residual allows to revisit parameter choice rules and relate a-priori, a-posteriori, and heuristic rules in a novel way that blurs the lines between the classical division of the parameter choice rules. All results are accompanied by numerical experiments.
\end{abstract}

\section{Introduction}\label{sec:problem}
In this paper we will revisit classical Tikhonov regularization. In this setting, we are interested in the solution of operator equations of the form
\begin{equation}\label{eq:problem}
y=Ax,
\end{equation}
where $A: X\rightarrow Y$ is a bounded, linear, and compact operator between Hilbert spaces $X$ and $Y$. Compactness implies that $A$ has a non-closed range, $\cR(A)\neq\overline{\cR(A)}$, and \eqref{eq:problem} is ill-posed. While there are operators with non-closed range that are not compact (strictly singular operators), we confine ourselves here to compactness since it allows to use the singular system of $A$ for analysis, and constitutes a natural limit for finite dimensional (and thus necessarily compact) approximations to $A$ used in numerical computation. For notational convenience we will consider $A$ to be scaled such that $\|A\|=1$. Instead of the exact data $y$ we have access only to a noisy datum $y^\delta$, for which we use the additive noise model 
\begin{equation}\label{eq:noise}
y^\delta=y+\epsilon \mbox{ with } \|\epsilon\|=\|y-y^\delta\|=\delta
\end{equation}
for a (typically unknown) noise level $\delta>0$. Tikhonov regularization approximates the unknown solution $\xd$ to \eqref{eq:problem} by solving
\begin{equation}\label{eq:tikh_classic_noise}
\xad=\arg\,\min \left\lbrace \frac{1}{2}\|Ax-y^\delta\|^2+\frac{\alpha}{2}\|x\|^2\right\rbrace,
\end{equation}
for $\alpha>0$. The optimal solution is known to be given as $\xad=(\asa+\alpha I)^{-1}A^\ast y^\delta$. Later we will sometimes consider noise-free data. In this case, we drop the superscript and denote the approximate solutions by $\xa$. Since the minimization problem is easily solved, the main task is to find an appropriate value of the regularization parameter $\alpha>0$. If $\alpha$ is too small, $\xad$ will follow the noise, whereas for $\alpha$ too large the solutions will be too smooth and thus also too far from $\xd$. Naturally, one would like to find the best possible approximation of $\xd$, i.e., minimize the reconstruction error $\|\xad-\xd\|$. This is quantified by convergence rates, which is a term for estimates of the form
\begin{equation}\label{eq:rates_prototype}
\|\xad-\xd\|\leq \varphi(\delta), \qquad 0<\delta\leq \delta_0
\end{equation}
with some \textit{index function} $\varphi$, i.e., $\varphi:[0,\infty)\rightarrow \R_+$ is continuous and monotonically increasing with $\varphi(0)=0$. With no further restriction on $\xd$, no such $\varphi$ exists \cite{EHN}. A classical assumption on $\xd$ is a \textit{source condition}, postulating the existence of a parameter $\mu>0$ such that 
\begin{equation}\label{eq:sc}
x^\dag\in\cR((A^\ast A)^\mu).
\end{equation}
Using the source condition, one can show that the optimal choice for the regularization parameter, given $\mu$ and $\delta$, is  
\begin{equation}\label{eq:apriori}
\alpha=c\delta^{\frac{2}{2\mu+1}}
\end{equation}
with a suitable constant $c>0$, and yields the convergence rate
\begin{equation}\label{eq:opt_rate_tikh}
\|\xad-\xd\|\leq C\delta^{\frac{2\mu}{2\mu+1}}
\end{equation}
for $0<\mu\leq 1$, and the exponent can not be reduced further. For more details on the convergence theory we refer to \cite{EHN}. A common remark on the theory sketched above is that the underlying assumptions are often difficult to verify in practice. The parameter choice \eqref{eq:apriori} and the estimate of the reconstruction error \eqref{eq:opt_rate_tikh} require the values of the smoothness parameter $\mu$ from \eqref{eq:sc} and the noise level $\delta$ \eqref{eq:noise}. While there are statistical methods that can potentially estimate the noise level, the source condition requires the knowledge of $\xd$, which is unavailable.

The main result of the paper is that, in the absence of noise, the asymptotics of the residual of the Tikhonov-regularized approximation to $\xd$ are equivalent to the source smoothness, namely, $\|A(\xa-\xd)\|=\bigo(\alpha^{\mu+\frac{1}{2}})$ for $0<\mu<\frac{1}{2}$ as $\alpha\rightarrow 0$ if and only if a variant of the source condition holds. By calculating approximations $\xad$ for various $\alpha$, we can trace the residual curve and, provided $\xd$ fulfils a source condition with $\mu<\frac{1}{2}$, extract the smoothness parameter by regression. We further show that once $\alpha$ is small enough such that the residual reaches the noise level, i.e., $\|A\xad-y^\delta\|\approx \delta$, decreasing $\alpha$ further will not significantly change the residual, which in turn allows inferring the noise level $\delta$ from the residual curve. Both parameter estimations, for $\mu$ and $\delta$, can be carried out for a single given datum $y^\delta$, making it applicable for any practical measurement.

A method for determining the smoothness parameter $\mu$ was first demonstrated in \cite{DGLoja}, based on exploiting a Kurdyka-\L{}ojasiewicz inequality implied by a source condition. However, the algorithm seemed instable, and some numerical observations remained unexplained. By using Tikhonov regularization instead of the Landweber method, we can calculate approximate solutions $\xad$ and the corresponding residuals for any regularization parameter $\alpha>0$, instead of being restricted to the discrete iteration steps, which increases the accuracy significantly.

As a byproduct of the noise level estimation we find a novel parameter choice rule for the regularization parameter $\alpha$, which we compare with other established parameter choice rules. Using our main result on the connection between residual asymptotics and solution smoothness, we can shed new light on the relation of parameter choice rules. In particular, we show that if $\xd$ satisfies a source condition with $0<\mu<\frac{1}{2}$, the a-priori parameter choice \eqref{eq:apriori}, the discrepancy principle as an a-posteriori parameter choice rule, and two heuristic rules, the heuristic discrepancy principle and our new method, differ only in the constant. Hence, all four yield order-optimal convergence rates, and we can blur the line between the three categories of parameter choice rules.

Since the source condition is difficult to verify in practice, we set up a simple model problem, for which we know the source condition of $\xd$. We use this model problem to illustrate our results.

\textbf{Model Problem:} Let $X=Y=\ell^2$, the space of square-summable sequences. We consider, for $\beta>0$ the operator $A:\ell^2\rightarrow \ell^2$, $[Ax]_i=i^{-\beta} x_i$, $i\in \N$. This is a compact operator with $\sigma_i=i^{-\beta}$, $i\in \N$, where the singular functions $v_i,u_i$ are the unit vectors in $\ell^2$. We set our exact solution $x^\dag=\{i^{-\eta}\}_{i=1}^\infty$ for some $\eta>\frac{1}{2}$. This yields $y=A\xd=\{i^{-\eta-\beta}\}_{i=1}^\infty$. Then, with $\mu^\ast=\frac{2\eta-1}{4\beta}$,
\begin{equation}\label{eq:sc_model}
x^\dag\in \bigcap_{\kappa<\mu^\ast} \cR((A^\ast A)^\kappa.
\end{equation}

The remainder of the paper is structured as follows. Our main result is contained in Section \ref{sec:conv}, where we show a converse result connecting source condition and the residual. We proceed by studying the noise in the residual in Section \ref{sec:noise}. The method for the estimation of solution smoothness and noise is summarized in Section \ref{sec:summary}. The converse result holds only for H\"{o}lder-type source conditions. In Section \ref{sec:lowhighorder} we discuss the case of higher and lower solution smoothness, which we show to be detectable in principle. A case study for real sets of tomographic data is presented in Section \ref{sec:tomo}. Finally, we discuss parameter choice rules in Section \ref{sec:paramchoice}.

\section{Converse Results}\label{sec:conv}
In \cite{Neubauer}, Neubauer showed that solution smoothness not only implies a convergence rate for the reconstruction error, but that also the reverse implication, often called a converse result, holds for Tikhonov regularization, both with noise-free and noisy data. The results were later generalized to other regularization approaches, see for example a generalization to Hilbert spaces \cite{AEHS2016,FHM2011} or Banach spaces \cite{F2018,MilHoh2020}, that all showed the equivalence of solution smoothness and convergence rates for the reconstruction error in their respective settings. We pursue here a different type of generalization. Instead of the relation between  reconstruction error and solution smoothness alone, we consider the behaviour of $\|\basa^\nu(\xa-\xd)\|$ as $\alpha\rightarrow 0$ and its relation to solution smoothness. Instead of the formulation \eqref{eq:sc} for the source condition, we follow Neubauer \cite{Neubauer} and use instead the condition
\begin{equation}\label{eq:sc_sum}
\sum_{n=k}^\infty \langle \xd,v_n\rangle^2=\bigo(\sigma_k^{4\mu}),
\end{equation}
for $k\rightarrow \infty$ which implies
\begin{equation}\label{eq:sc_open}
x^\dag\in \bigcap_{\kappa<\mu} \cR((A^\ast A)^\kappa).
\end{equation}
Note that this is the setting of our model problem, see \eqref{eq:sc_model}.

Since $A$ is assumed to be compact, we can use its singular system $\{\sigma_i,u_i,v_i\}_{i=1}^\infty$ for the analysis. There, the functions $\{u_i\}_{i=1}^\infty$ form an ONB for $\overline{\cR(A)}$, $\{v_i\}_{i=1}^\infty$ form an ONB for $\overline{\cR(A^\ast)}$ and the singular values $\{\sigma_i\}_{i=1}^\infty$ accumulate at zero; provided $\mathrm{dim}(\cR(A))=\infty$. We recall that the relations $Av_i=\sigma_i u_i$ and $A^\ast u_i=\sigma_i v_i$ hold for all $i\in \N$. Any $x\in X$ can be written as $x=\sum_{i=1}^\infty \langle x,v_i\rangle v_i$ and 
\[
Ax= \sum_{i=1}^\infty \langle Ax,u_i\rangle u_i=\sum_{i=1}^\infty \langle x,v_i\rangle Av_i .
\]

The following lemma is the basis for our converse result.
\begin{lemma}\label{lem:sumeval}
Let $\{\sigma_i,v_i,u_i\}_{i=1}^\infty$ be the singular system to $A$, $\xd\in X$. Then 
\begin{equation}
\sum_{i=1}^\infty \frac{\sigma_i^{q}\lambda^2}{(\sigma_i^{p}+\lambda)^2}\langle \xd,v_i\rangle^2 = \bigo(\lambda^{\frac{q+4\mu}{p}})
\end{equation}
for $0<q+4\mu<2p$ if and only if $\xd$ satisfies \eqref{eq:sc_sum}.
\end{lemma}
\begin{proof}
In principle the proof follows that of \cite[Theorem 1]{Neubauer}. We split
\begin{align}\label{eq:lemma_split}
\sum_{i=1}^\infty \frac{\sigma_i^{q}\lambda^2}{(\sigma_i^{p}+\lambda)^2}\langle \xd,v_i\rangle^2=\sum_{\sigma_i^p\leq\lambda}\frac{\sigma_i^{q}\lambda^2}{(\sigma_i^{p}+\lambda)^2}\langle \xd,v_i\rangle^2+\sum_{\sigma_i^p> \lambda}\frac{\sigma_i^{q}\lambda^2}{(\sigma_i^{p}+\lambda)^2}\langle \xd,v_i\rangle^2.
\end{align}
Consider first the small singular values. Noting that $\frac{1}{4}\leq \frac{\lambda^2}{(\sigma_i^p+\lambda)^2}\leq 1$ for $\sigma_i^p<\lambda$, we have
\[
\sum_{\sigma_i^p\leq\lambda}\frac{\sigma_i^{q}\lambda^2}{(\sigma_i^{p}+\lambda)^2}\langle \xd,v_i\rangle^2\leq \lambda^{\frac{q}{p}}\sum_{\sigma_i^p\leq\lambda} \langle \xd,v_i\rangle^2=\bigo(\lambda^{\frac{q}{p}}(\lambda^{\frac{1}{p}})^{4\mu})=\bigo(\lambda^{\frac{q+4\mu}{p}}).
\]

For the term corresponding to the larger singular values, we use that $\frac{1}{4}\leq\frac{\sigma^{2p}}{(\sigma^p+\lambda)^2}<1$ for $\lambda\leq \sigma\leq 1$, which yields
\[
\frac{\lambda^2}{4}\sum_{\sigma_i^p>\lambda}\sigma_i^{q-2p}\langle \xd,v_i\rangle^2\leq \sum_{\sigma_i^p>\lambda} \sigma_i^{q-2p}\lambda^2\frac{\sigma_i^{2p}}{(\sigma_i^p+\lambda)^2}\langle \xd,v_i\rangle^2< \lambda^2\sum_{\sigma_i^p>\lambda}\sigma_i^{q-2p}\langle \xd,v_i\rangle^2.
\]
It remains to show $\lambda^2\sum_{\sigma_i^p>\lambda}\sigma_i^{q-2p}\langle \xd,v_i\rangle^2=\bigo(\lambda^{\frac{q+4\mu}{p}-2})$. Via induction (see Appendix) one finds with \eqref{eq:sc_open} that
\[
\sum_{i=1}^k \sigma_i^{q-2p}\langle\xd,v_i\rangle^2=\bigo(\sigma_k^{q+4\mu-2p}).
\]
Now we choose $k$ such that $\sigma_k^p=\bigo(\lambda)$. With $q+4\mu-2p\leq0$ it follows
\[
\lambda^2\sum_{\sigma_i^p>\lambda}\sigma_i^{q-2p}\langle \xd,v_i\rangle^2=\lambda^2\bigo(\sigma_k^{q+4\mu-2p})=\lambda^2\bigo(\lambda^{\frac{q+4\mu}{p}-2})=\bigo(\lambda^{\frac{q+4\mu}{p}}).
\]
This, together with the upper bound for the first summand in \eqref{eq:lemma_split}, yields the claim.
\end{proof}

\begin{theorem}\label{thm:conv}
Let $0<\mu+\nu<1$. Then
\begin{equation}\label{eq:numurates}
\|\basa^\nu(\xa-\xd)\|^2=\bigo\left(\alpha^{2(\nu+\mu)}\right)
\end{equation}
if and only if  $\xd$ satisfies \eqref{eq:sc_sum}.
\end{theorem}
\begin{proof}
It is
\[
\|\basa^\nu(\xa-\xd)\|^2=\sum_{\sigma_i} \frac{\sigma_i^{4\nu}\alpha^2}{(\sigma_i^2+\alpha)^2}\langle \xd,v_i\rangle^2.
\]
Hence we apply Lemma \ref{lem:sumeval} with $q=4\nu$ and $p=2$.

\end{proof}

Neubauer \cite{Neubauer}, as well as Scherzer et. al. \cite{AEHS2016} also provide converse results for noisy data. Since we assume that we only have one set of data $y^\delta$ with fixed $\delta$ available, we will not pursue convergence rates for noisy data further. Due to the additive noise model and Theorem \ref{thm:conv}, the equivalence between \eqref{eq:sc_sum} and order optimal convergence rate would be no surprise.

Theorem \ref{thm:conv} states that the solution smoothness in Tikhonov regularization is preserved under the application of certain powers of $\asa$. Hence, the smoothness of the powers $\basa^\nu(\xa-\xd)$ can be used to assess the solution smoothness. In practice, most values of $\nu$ are still not observable. However, for $\nu=\frac{1}{2}$ we find the residual $\|\basa^\frac{1}{2}(\xa-\xd)\|=\|A(\xa-\xd)\|$, and for $\nu=1$ we obtain the gradient. For Tikhonov regularization, all solution smoothness is lost in the gradient, since $\nu+\mu<1$ is a requirement of Theorem \ref{thm:conv}. The loss of information in the gradient can also be seen from the first order condition
\begin{equation}\label{eq:firstorder}
\asa(\xa-\xd)=-\alpha\xa,
\end{equation}
which enforces that $\|\asa(\xa-\xd)\|= \alpha\|\xa\|$ unconditionally. Because $\xa\rightarrow \xd$ as $\alpha\rightarrow 0$ (note that we are in the noise-free scenario), $\|\xa\|\rightarrow \|\xd\|$ and hence $\|\asa(\xa-\xd)\|\sim \alpha$.

The residual, obtained with $\nu=\frac{1}{2}$, on the other hand, contains smoothness information. Since Theorem \ref{thm:conv} requires $\mu+\nu<1$, smoothness with $\mu<\frac{1}{2}$ is preserved and hence can be detected. It is no coincidence that this matches with the well-known fact that the discrepancy principle for the choice of $\alpha$ yields order-optimal convergence rates for $0<\mu<\frac{1}{2}$, as discussed in more detail in a follow-up paper.

Another interesting observation is that the saturation of Tikhonov regularization, i.e., the fact that the best obtainable convergence rate is $\|\xad-\xd\|=\bigo(\delta^{\frac{2}{3}})$ for $\mu\geq1$ is due to the effect that low-frequency components of the solutions cannot be approximated well: The saturation follows from the condition $\mu+\nu<1$ in Theorem \ref{thm:conv}, which we have used solely to evaluate the low frequencies $\sigma_i^2>\alpha$.

\section{Data noise}\label{sec:noise}
\begin{figure}
\includegraphics[width=\linewidth]{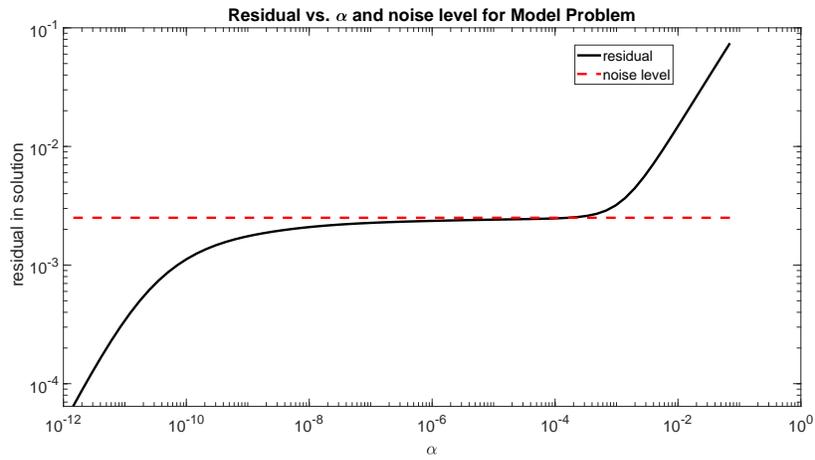}\caption{Residual vs. $\alpha$ and noise level $\delta$ for Tikhonov regularization of our Model Problem with $\eta=\beta=2$ ($\mu=0.375$). For large $\alpha$, the approximation dominates, hence the residual decreases with $\alpha$. Even more, in this phase there holds $\|A\xad-y^\delta\|\sim \|A\xa-y\|\sim \alpha^{\mu+0.5}$, i.e, the residual follows the the theoretical behaviour for noise-free data. Once the residual is close to the noise level, it begins to stagnate around the noise level: 5 magnitudes of regularization parameters yield an almost identical residual value. Only for the smallest $\alpha$ the residuals decay strongly again.}\label{fig:resexample}
\end{figure}

So far we have not considered noise in the data. Traditionally, the analysis is focused on the propagation of the noise to the reconstruction error. Instead, here we focus again on the residuals. Recall that we assume the additive noise model \eqref{eq:noise}.
Using the singular system we express the residual as
\begin{align}
A\xad-y^\delta&=\sum_{i=1}^\infty \left(\frac{\sigma_i^2}{\sigma_i^2+\alpha}-1\right)\langle y^\delta,u_i\rangle u_i\nonumber\\
&=\sum_{i=1}^\infty \frac{\alpha}{\sigma_i^2+\alpha}\langle y,u_i\rangle u_i+\sum_{i=1}^\infty \frac{\alpha}{\sigma_i^2+\alpha}\langle \epsilon,u_i\rangle u_i\label{eq:res_1}\\
&=\sum_{i=1}^\infty \frac{\sigma_i\alpha}{\sigma_i^2+\alpha}\langle \xd,v_i\rangle u_i+\sum_{i=1}^\infty \frac{\alpha}{\sigma_i^2+\alpha}\langle \epsilon,u_i\rangle u_i.\label{eq:res_2}
\end{align}
Let us first look at \eqref{eq:res_1}. The factors for the data approximation (left term) and noise (right term) are identical. Since $y$ and $\epsilon$ are fixed, this means that, without further information, the filtering of the exact data and of the noise behave identically when $\alpha$ is varied. Therefore it is, in the most general setting of arbitrary $y,\epsilon\in Y$, not possible to distinguish between the exact data $y=A\xd$ and the noise $\epsilon$. In the following we explicitly exclude the theoretical case that the exact data $y^\delta$ and the noise have similar smoothness and assume that
\[
\langle y,u_i\rangle = o( \langle \epsilon,u_i\rangle).
\]
Since $y\in \cR(A)$ due to \eqref{eq:problem}, it possesses a minimal smoothness which is in most situations enough to distinguish it from the noise. In \eqref{eq:res_2} we have inserted \eqref{eq:problem} into \eqref{eq:res_1} and used that $A^\ast u_i=\sigma_iv_i$. Now the filters for approximation and noise are no longer identical. Let $\alpha$ be fixed. Then $\frac{\sigma_i\alpha}{\sigma_i^2+\alpha}\rightarrow 0$ as $\sigma_i\rightarrow 0$, but $\frac{\alpha}{\sigma_i^2+\alpha}\rightarrow 1$ at the same time. This means that the high-frequency components of the approximation are suppressed, while the high-frequency components of the noise remain almost unchanged. Assuming that $\sigma_i\langle x^\dag,v_i\rangle$ decays much faster than $\langle\epsilon,u_i\rangle$, the noise will be dominating the residual, and $\|A\xad-y^\delta\|\approx \delta$ for a large range of regularization parameters $\alpha$. In Figure \ref{fig:resexample} we demonstrate this by using our Model Problem contaminated by Gaussian noise. Note that in a discretized setting, i.e., when summing only up to some finite $N\in\N$, the noise will be suppressed for sufficiently small $\alpha$. Namely, due to the discretization, the sums in \eqref{eq:res_1} for the residual are truncated at some $N>0$. Because the factors $\frac{\alpha}{\sigma_i^2+\alpha}\rightarrow 0$ as $\alpha\rightarrow 0$, the residual eventually behaves as $\bigo(\alpha)$ for $\alpha$ very small. More precisely, since $\sigma_i\geq \sigma_N$ for $i=1,\dots,N$, 
\[
\sum_{i=1}^N \frac{\alpha}{\sigma_i^2+\alpha}\langle \epsilon,u_i\rangle \approx \alpha \sum_{i=1}^N \frac{\langle \epsilon,u_i\rangle}{\sigma_i^2}
\]
if $\alpha$ is sufficiently small. 

The observations above can be used to estimate the noise level by looking for the flat plateau in the residual. One possibility of doing this is to look for the saddle point in the residual curve, cf. Figure \ref{fig:resexample}. We make this more precise in the next section.

In a last remark on data noise we mention that in practical applications one will often encounter another form of noise: the modelling error. The practical measurement setups are to some extent idealized in the mathematical model. For example ray sources and detectors are modelled as points, although they have a small area in practice. Also their position cannot be measured to arbitrary precision. That means a practical measurement $y^{meas}$ will not coincide with $y_m^\delta$ as expected from the modelling. Even more, the case $y^{meas}\notin \overline{\cR(A)}$ is to be expected. The resulting difference $\delta^{meas}:=\|y^{meas}-y_m^\delta\|$ thus can not be explained by the model and remains as a limit of the residual, $\|A\xad-y^{meas}\|\geq \delta^{meas}$ for all $\alpha\rightarrow0$. In Section \ref{sec:tomo} we illustrate this numerically.

\section{Estimation of smoothness parameter and noise level}\label{sec:summary}
We now summarize our results to estimate solution smoothness and noise level using Tikhonov regularization.

The method for the estimation of the source condition is based on Theorem \ref{thm:conv}. Hence we can only identify the smoothness parameter if $\xd\in\cR(\basa^\nu)$ with $\nu<\frac{1}{2}$. The method works as follows. 
We employ Tikhonov regularization for several magnitudes of regularization parameters, and store the residuals 
\[
r(\alpha):=\|A\xad-y^\delta\|.
\]
We also compute the derivative
\[
dr(\alpha):=\frac{\partial}{\partial \log(\alpha)}\log r(\alpha).
\]

\begin{figure}
\includegraphics[width=\linewidth]{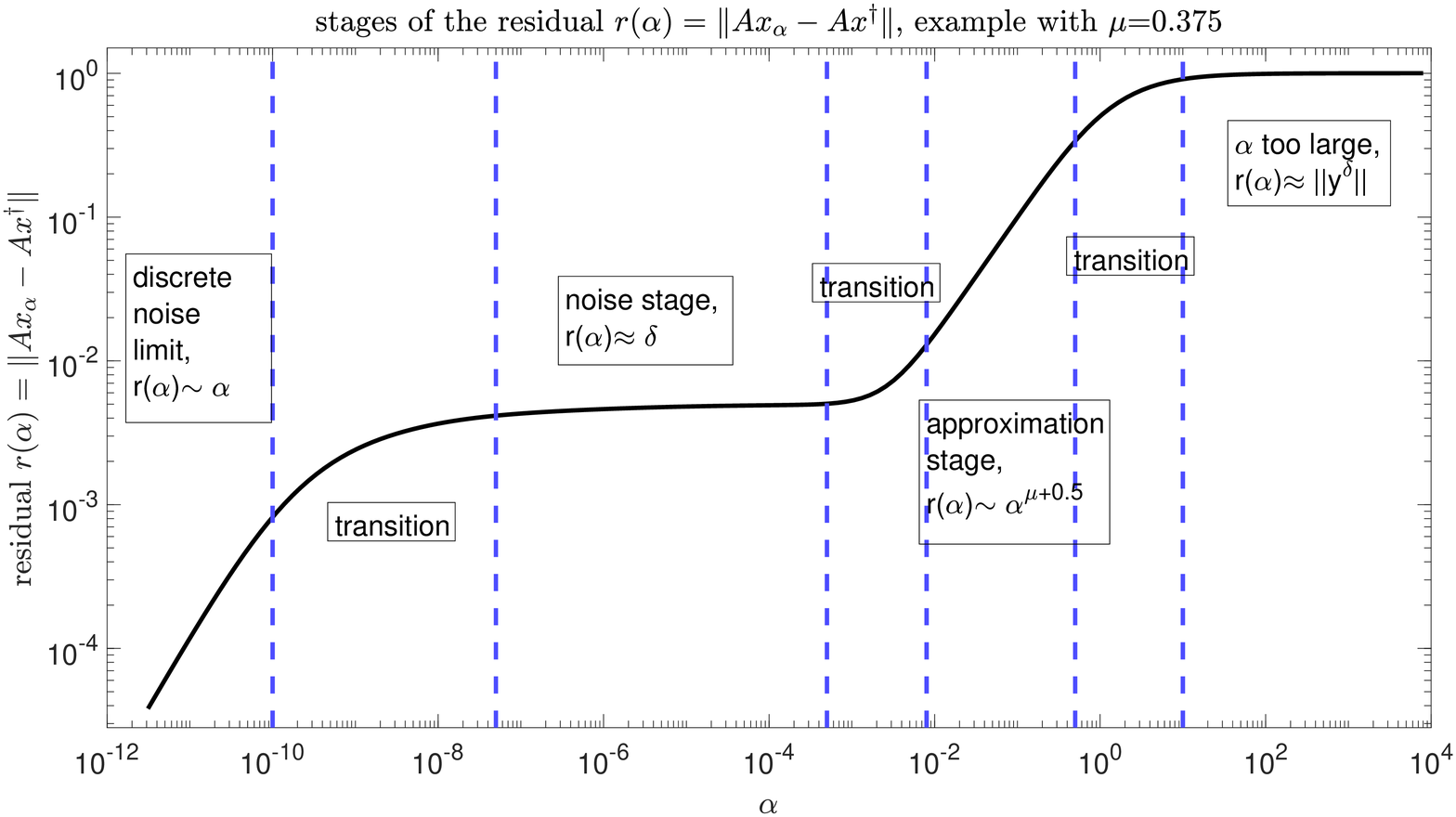}\\\includegraphics[width=\linewidth]{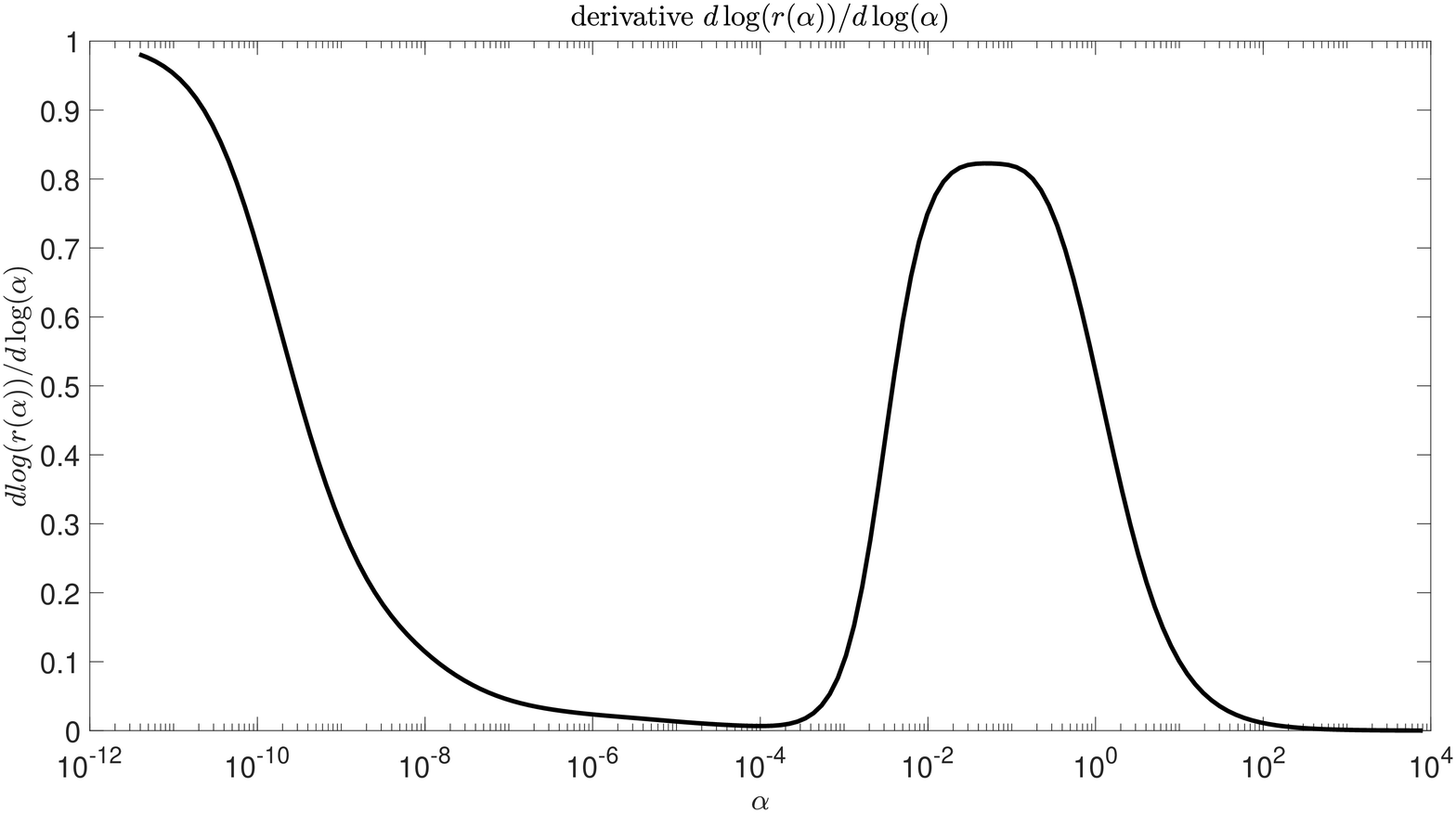}\caption{Top: Example of the stages of the residual under Tikhonov regularization for a large range of regularization parameters. Bottom: corresponding derivative of the residual.}\label{fig:residualcurve}
\end{figure}

Both curves as function of $\alpha$ can be characterized in different intervals or stages, that can be used to estimate the smoothness parameter $\mu$ and the point to extract an estimate of the noise level. An example for the Model Problem with $\eta=\beta=2$ and $\delta=0.005\|y\|$ is given in Figure \ref{fig:residualcurve}. 

When the regularization parameters are too large, we have $\|\xad\|\approx 0$, hence $r(\alpha)\approx \|y^\delta\|$ and $dr(\alpha)\approx0$. Lowering $\alpha$ and going through a transition stage, we arrive at the approximation phase. Due to \eqref{eq:res_1} it is $r(\alpha)\leq\|A\xa-y\|+\|y-y^\delta\|$, so, since $\|y-y^\delta\|=\delta$, if $r(\alpha)>>\delta$ this means $r(\alpha)\approx \|A\xa-y\|$. The residual is dominated by the approximation of the exact data and the noise has little impact on the residual, which therefore carries the information on solution smoothness. Theorem \ref{thm:conv} with $\nu=\frac{1}{2}$ yields $r(\alpha)\sim \alpha^{\mu+\frac{1}{2}}$ if and only if $\eqref{eq:sc_open}$. Hence, we can make a regression for the ansatz $r(\alpha)=c\alpha^\kappa$. If the regression yields a good fit to the residual curve with $\frac{1}{2}<\kappa<1$, we can extract the solution smoothness $\xd\in \cR(\basa^{\mu^\ast})$ with $\mu^\ast=\kappa-\frac{1}{2}$. Due to numerical inaccuracies care has to be taken when $\mu^\ast\approx 0$ or $\mu^\ast\approx \frac{1}{2}$, since one might mistake low or high H\"{o}lder smoothness for logarithmic smoothness or, respectively, smoothness higher than an Hölder source condition with $\mu=\frac{1}{2}$, see Section \ref{sec:lowhighorder}. In the derivative $dr(\alpha)$, we can immediately read off $\mu$: if $dr(\alpha)=\kappa$ for some $\frac{1}{2}<\kappa<1$, then $\xd$ must satisfy the source condition \eqref{eq:sc_sum} with $\mu=\kappa-\frac{1}{2}$.

After another transition, we are in the noise phase, where $r(\alpha)\approx\delta$. This is again a consequence of \eqref{eq:res_1}, $r(\alpha)\leq\|A\xa-y\|+\|y-y^\delta\|$, and the comments thereafter. For the middle term it still holds $\|A\xa-y\|=\bigo(\alpha^{\mu+\frac{1}{2}})$ if \eqref{eq:sc_sum} holds, whereas $\|y-y^\delta\|=\delta$. Hence, whenever $\|A\xa-y\|\leq \delta$, $\|y-y^\delta\|$ is dominating the residual and the approximation component $\|A\xa-y\|$ is neglectable; $\|A\xad-y^\delta\|\approx \delta$ for a large range of regularization parameters $\alpha$. The noise level can therefore be read off the flat part of the residual, cf. Figure \ref{fig:resexample} and Figure \ref{fig:residualcurve}. From the figures we also see that the residual curve has almost a saddle point in the flat plateau, which we can use to find it and estimate the noise level algorithmically. We look for
\begin{equation}\label{eq:minalpharescurve}
\alpha^\ast=\argmin dr(\alpha)=\argmin \frac{\partial}{\partial \log(\alpha)}\log(\|A\xad-y^\delta\|)
\end{equation} 
and estimate $\delta\approx \|Ax_{\alpha^\ast}^\delta-y^\delta\|$. We discuss \eqref{eq:minalpharescurve} as a rule to choose the regularization parameter $\alpha$ in Section \ref{sec:paramchoice}.

In theory, the noise phase is then active for all $\alpha\rightarrow 0$ and $\|A\xad-y^\delta\|$ slowly goes to zero, depending on the decay of the noise components $\langle y-y^\delta,u_i\rangle$. In practice, due to discretization, we observe one or two more stages as $\alpha$ is decreased further. As noted in Section \ref{sec:noise}, we can expect to see $r(\alpha)\approx \alpha$ and thus $dr(\alpha)\approx 1$ for sufficiently small $\alpha$. It might happen that, when reducing the regularization parameter even further, the residuals become chaotic, likely due to numerical errors such as round-off errors and the amplification thereof due to the ill-posed nature of the problem.

We finally note that in order for this estimation to work we require that the noise level be not too high, i.e., $\delta<<\|y^\delta\|$. If it is, then the approximation stage is too short or even non-existent, such that the solution smoothness is completely hidden in the noise.

The algorithm for the estimation of $\mu$ and $\delta$ is summarized in Algorithm \ref{alg}. Note that the closer $q$ is to one, the clearer the expected results are.

\begin{algorithm}\caption{Algorithm for the estimation of $\mu$ and $\delta$.}
\begin{algorithmic}
\STATE Input: $A$, $y^\delta$, $\alpha_0>0$ such that $r(\alpha_0)\approx\|y^\delta\|$, $1<q<1$, $n\in N$
\FOR{$i=1,2,\dots,N$}
\STATE set $\alpha_i=\alpha_0 q^{i-1}$
\STATE calculate $x_{\alpha_i}^\delta=(A^\ast A+\alpha_i I)^{-1}A^\ast y^\delta$
\STATE store $r(\alpha_i)=\|x_{\alpha_i}^\delta-y^\delta\|$
\STATE calculate $dr(\alpha_i)\approx\frac{\log(r(\alpha_i))-\log(r(\alpha_{i-1}))}{\log(\alpha_i)-\log(\alpha_{i-1})}$
\ENDFOR
\STATE If $dr(\alpha_i)\approx \kappa$ with $\frac{1}{2}<\kappa<\mu$ for sufficiently many $\alpha_i$, estimate $\mu\approx \kappa-\frac{1}{2}$.
\STATE find $\alpha^\ast=\argmin dr(\alpha_i)$, estimate $\delta\approx \|Ax_{\alpha^\ast}^\delta-y^\delta\|$
\end{algorithmic}\label{alg}
\end{algorithm}

\section{Low smoothness and high smoothness}\label{sec:lowhighorder}
The estimation of the solution smoothness works best when the classical H\"{o}lder-type source conditions \eqref{eq:sc} or \eqref{eq:sc_open} describe the smoothness of $\xd$, as Theorem \ref{thm:conv} can be applied. More general, for each $\xd\in X$ there exists an index function $\varphi$ and $w\in X$ such that
\begin{equation}\label{eq:sc_general}
\xd=\varphi(\asa) w,
\end{equation}
see \cite{HofMathprofile}. We distinguish two cases, depending on whether $\varphi(t)$ decays slower or faster to zero than the H\"{o}lder-type functions $t^\kappa$, $0<\kappa<\frac{1}{2}$. In analogy to \eqref{eq:sc_sum}, we consider the generalized source conditions of the form 
\begin{equation}\label{eq:sc_general_sum}
\sum_{n=k}^\infty \langle x^\dag,v_n\rangle^2=\bigo(\varphi(\sigma_k^2)^2).
\end{equation}

\begin{figure}
\includegraphics[width=0.49\linewidth]{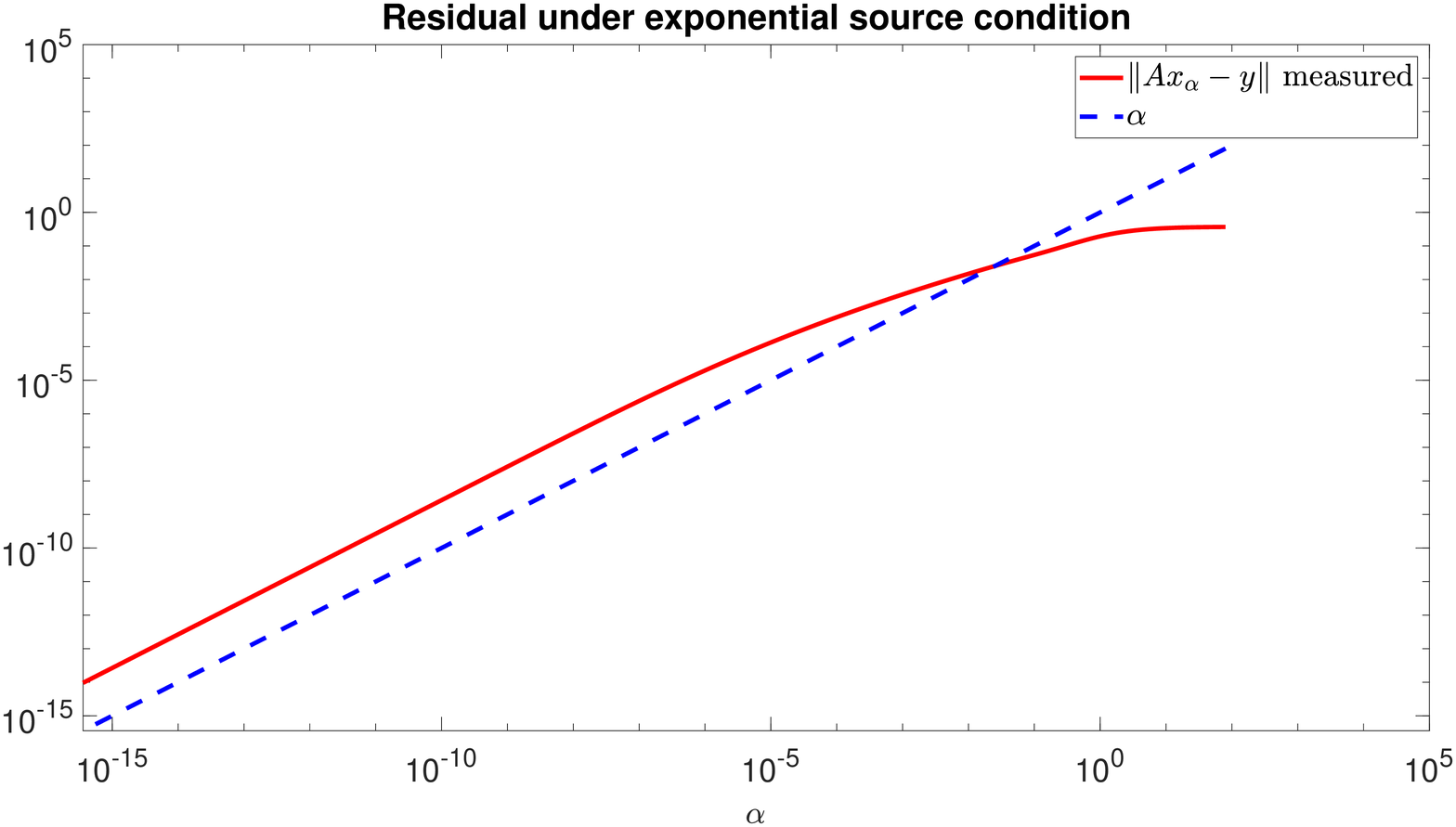}\includegraphics[width=0.49\linewidth]{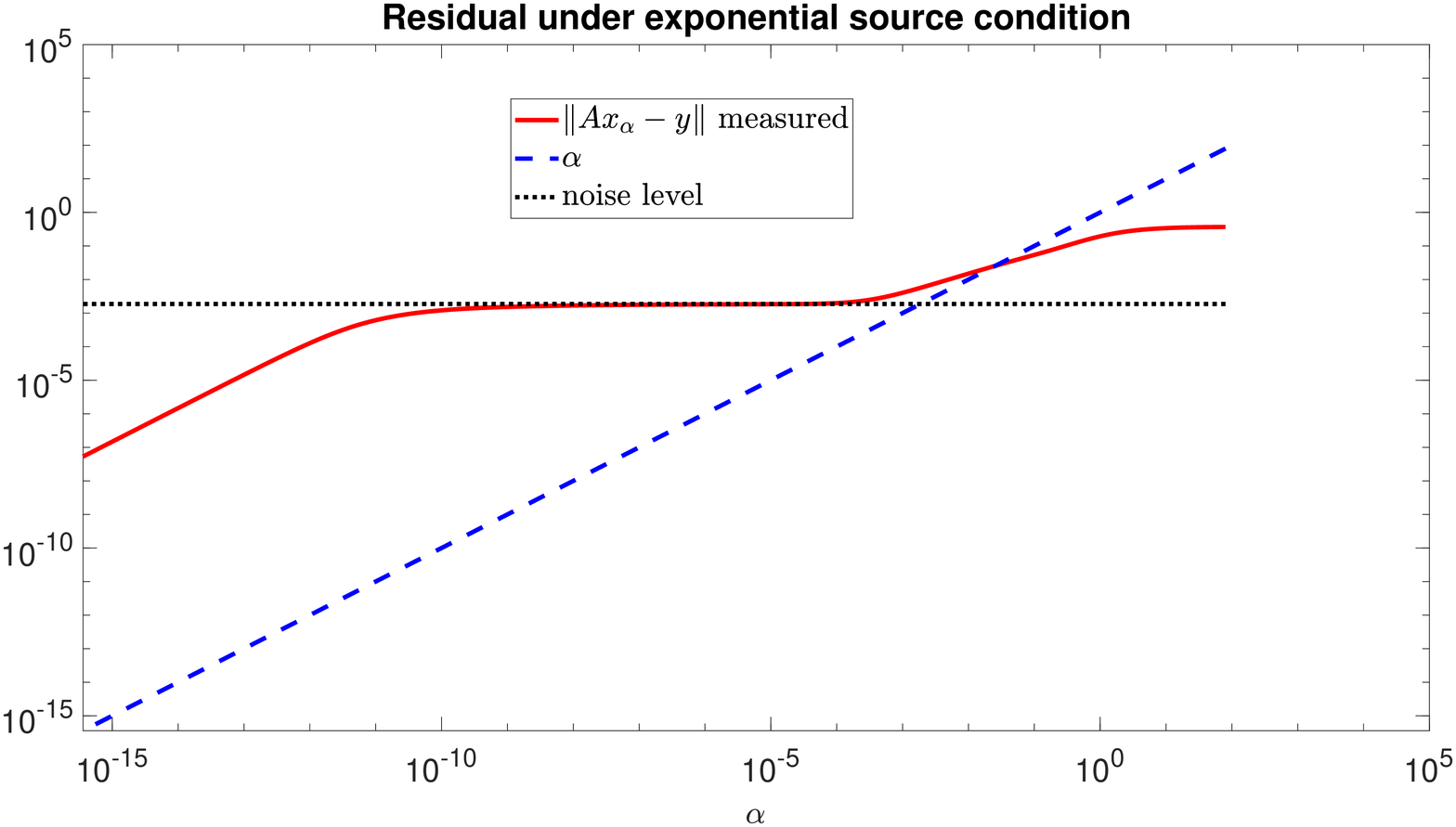}\\
\includegraphics[width=0.49\linewidth]{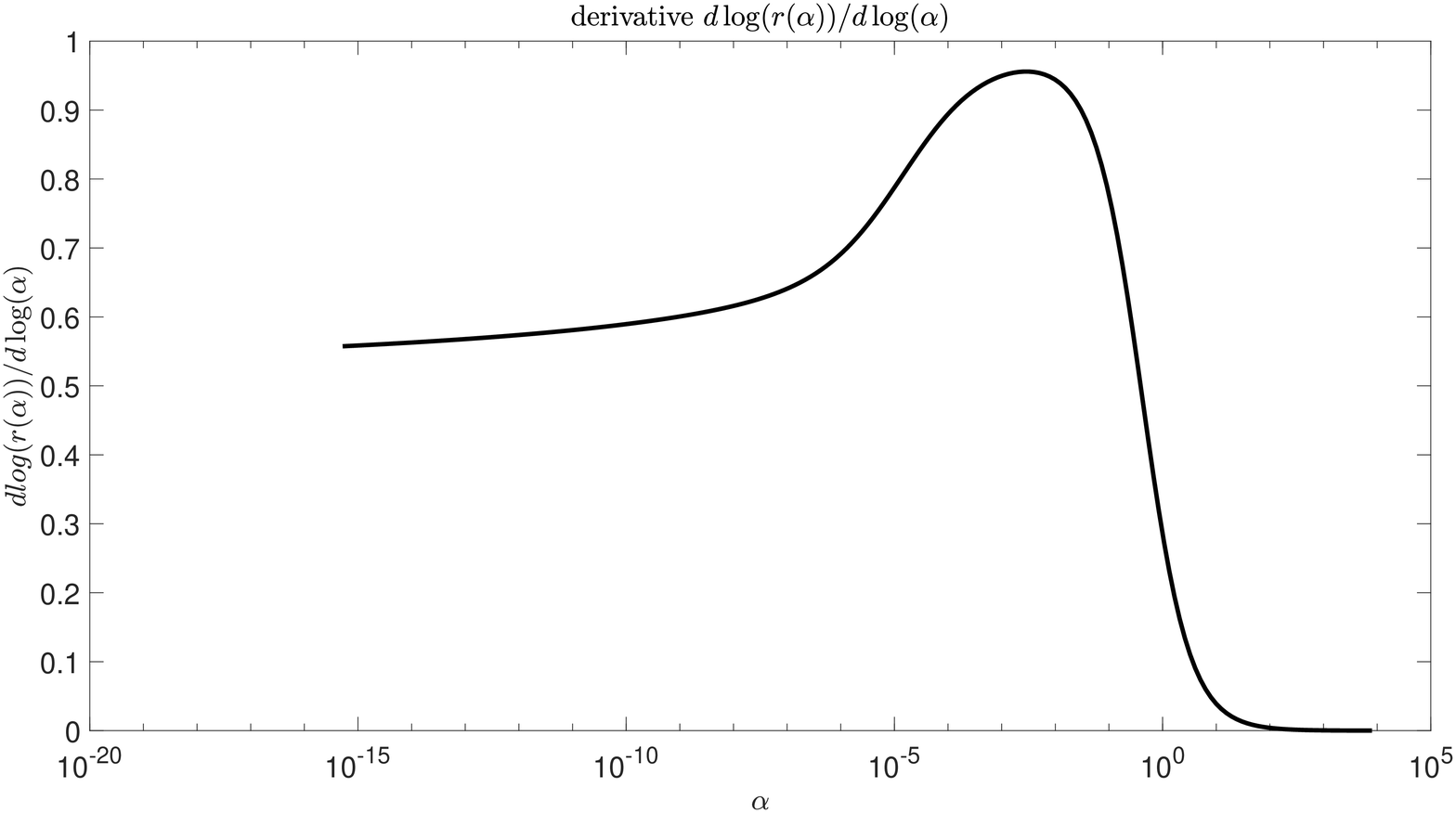}\includegraphics[width=0.49\linewidth]{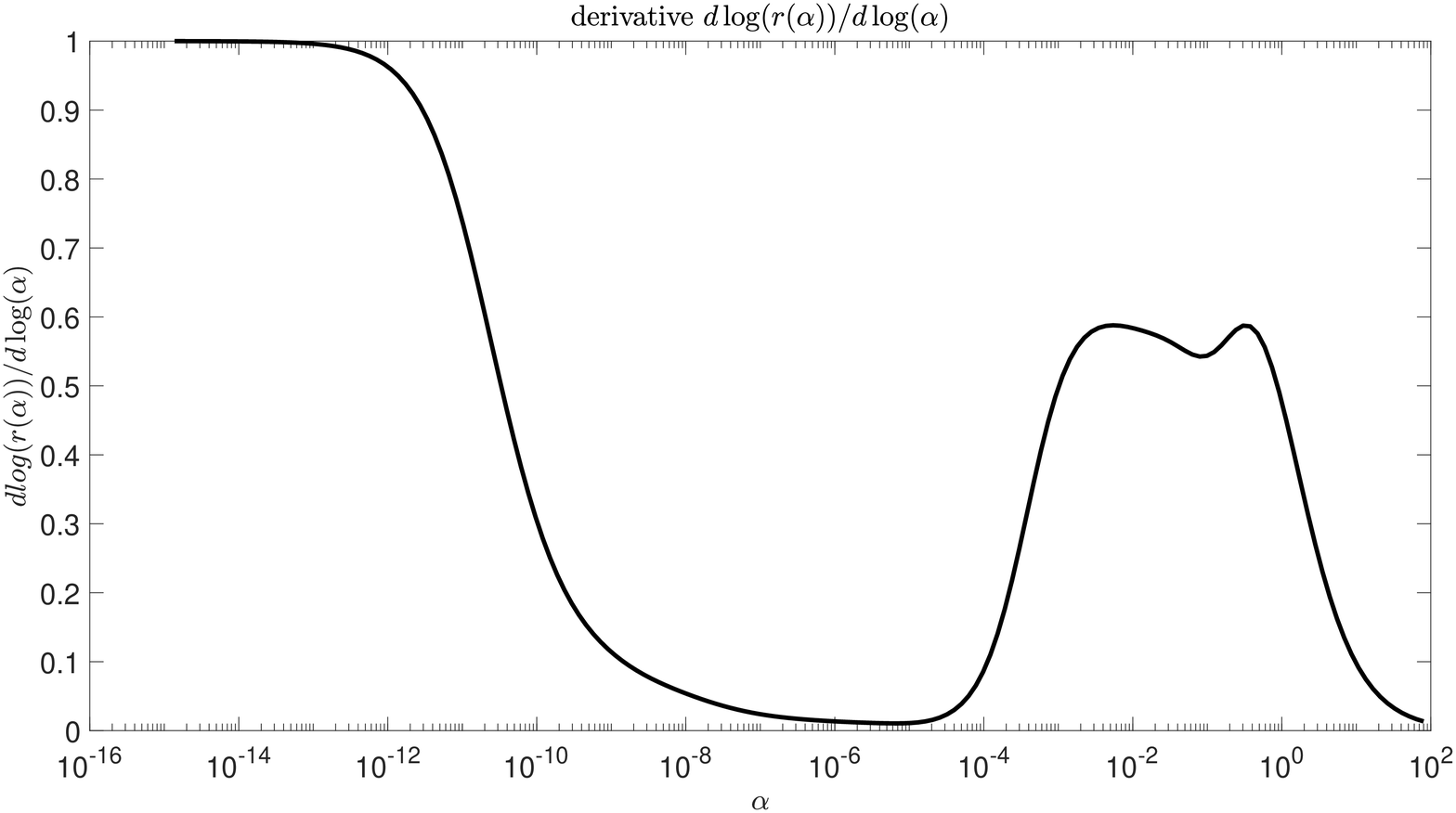}\caption{Numerical Results for the exponential smoothness setting \eqref{eq:sc_general_sum} with $\varphi(t)=\exp(-t^\frac{1}{2})$. Top: residuals, bottom: corresponding derivative. Because $\xd$ is smoother than a H\"{o}lder source condition with $\mu=\frac{1}{2}$, the smoothness is eventually lost in the residual and $\|A\xa-A\xd\|\sim \alpha$ for small enough $\alpha$ in the noise free case (left). For noisy data (right) we see the plateau indicating the noise level. In both cases we see that for the larger $\alpha$ the residual curve is concave.}\label{fig:expsmooth}
\end{figure}

We have seen in Theorem \ref{thm:conv} that the residual carries no solution smoothness information when $\xd$ fulfils a source condition \eqref{eq:sc_open} with $\mu\geq\frac{1}{2}$. This remains the case for functions smoother than the H\"{o}lder powers, as then asymptotically $\|A\xa-A\xd\|=\bigo( \alpha)$. A numerical example with an exponential source condition $\varphi(t)= \exp(-t^\frac{1}{\kappa})$ with $\kappa=2$ is shown in Figure \ref{fig:expsmooth}. In absence of noise the residual curve is, for larger regularization parameters, concave in the $\log$-$\log$ plot but for small enough $\alpha$ it is of order $\alpha$ as expected. For noisy data, the latter phase is completely masked by the noise. Still, one can easily spot the noise level. The visible part of the residual curve for larger $\alpha$ (i.e., the approximation phase) corresponds to the concave part of the noise-free residual, which appears to be the indicator for the high smoothness case.

\begin{figure}
\includegraphics[width=0.49\linewidth]{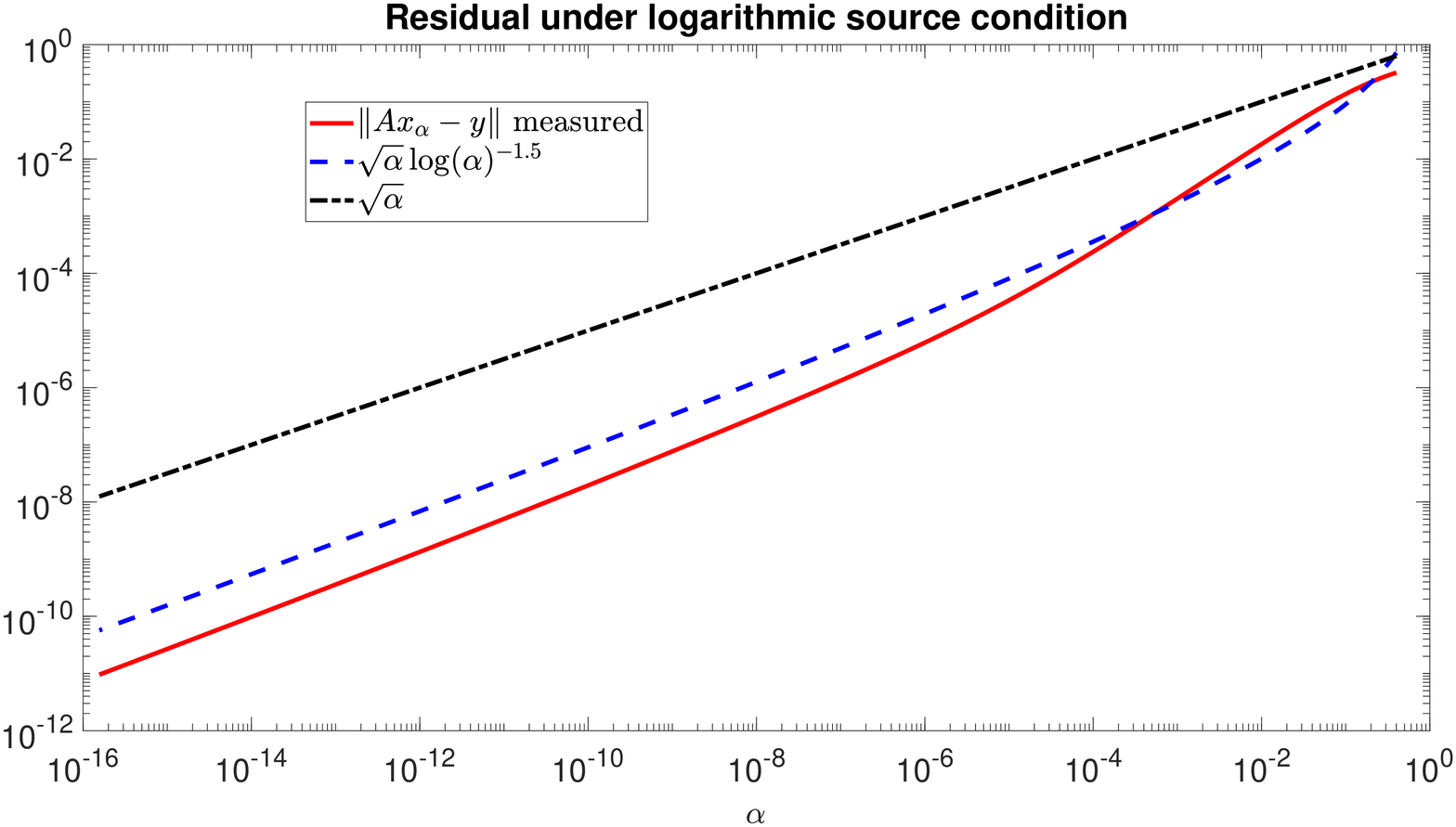}\includegraphics[width=0.49\linewidth]{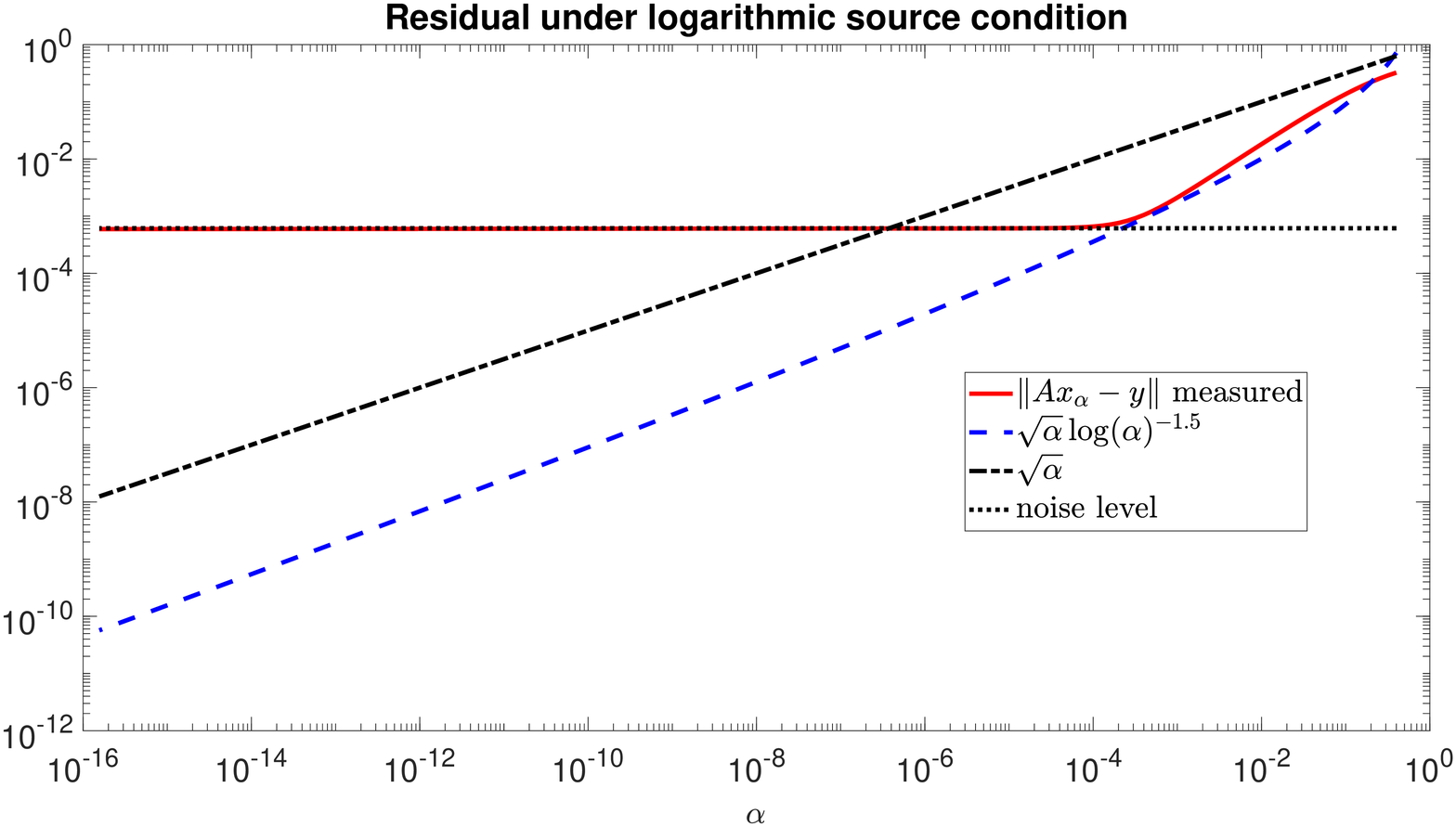}\\
\includegraphics[width=0.49\linewidth]{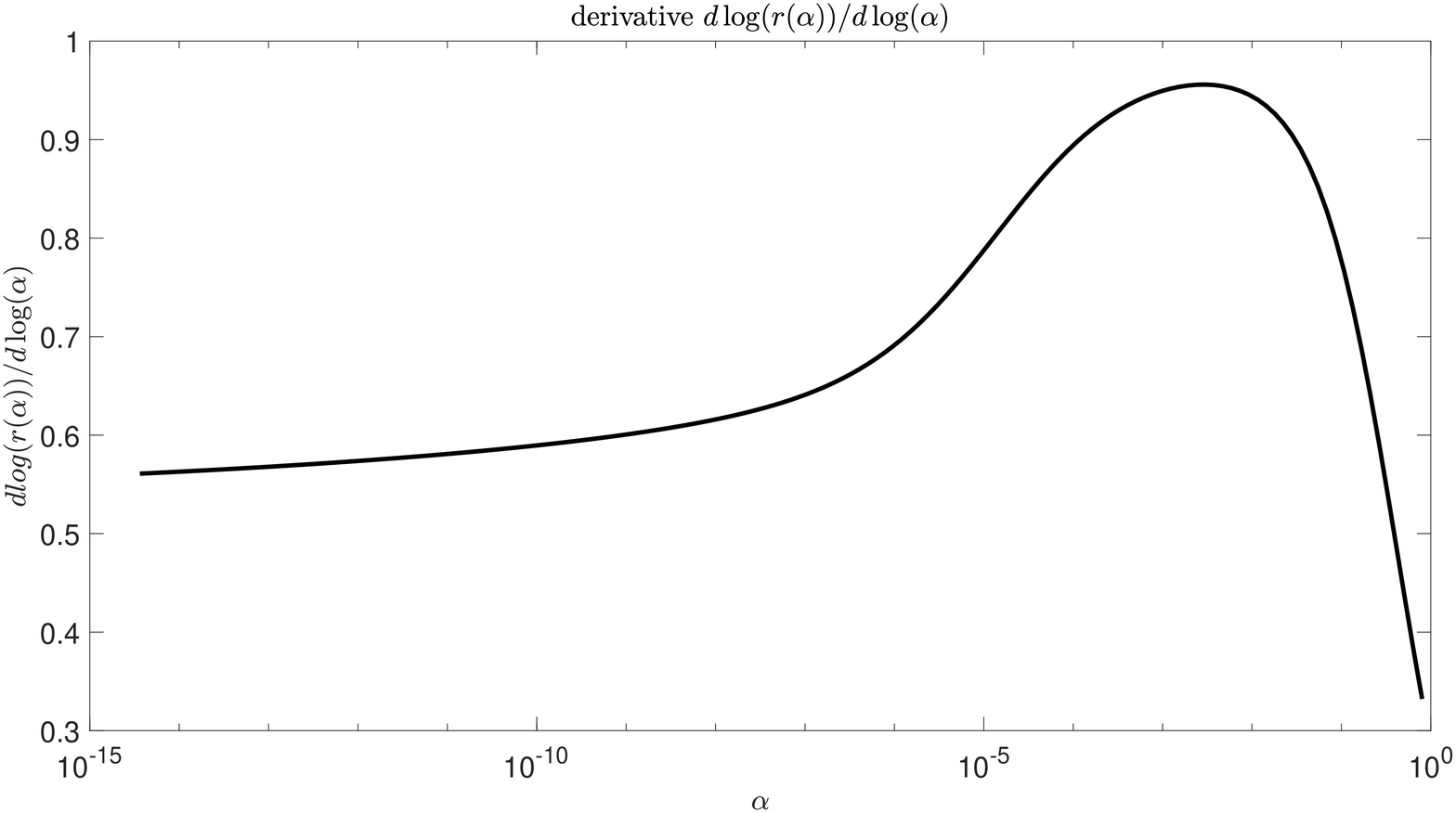}\includegraphics[width=0.49\linewidth]{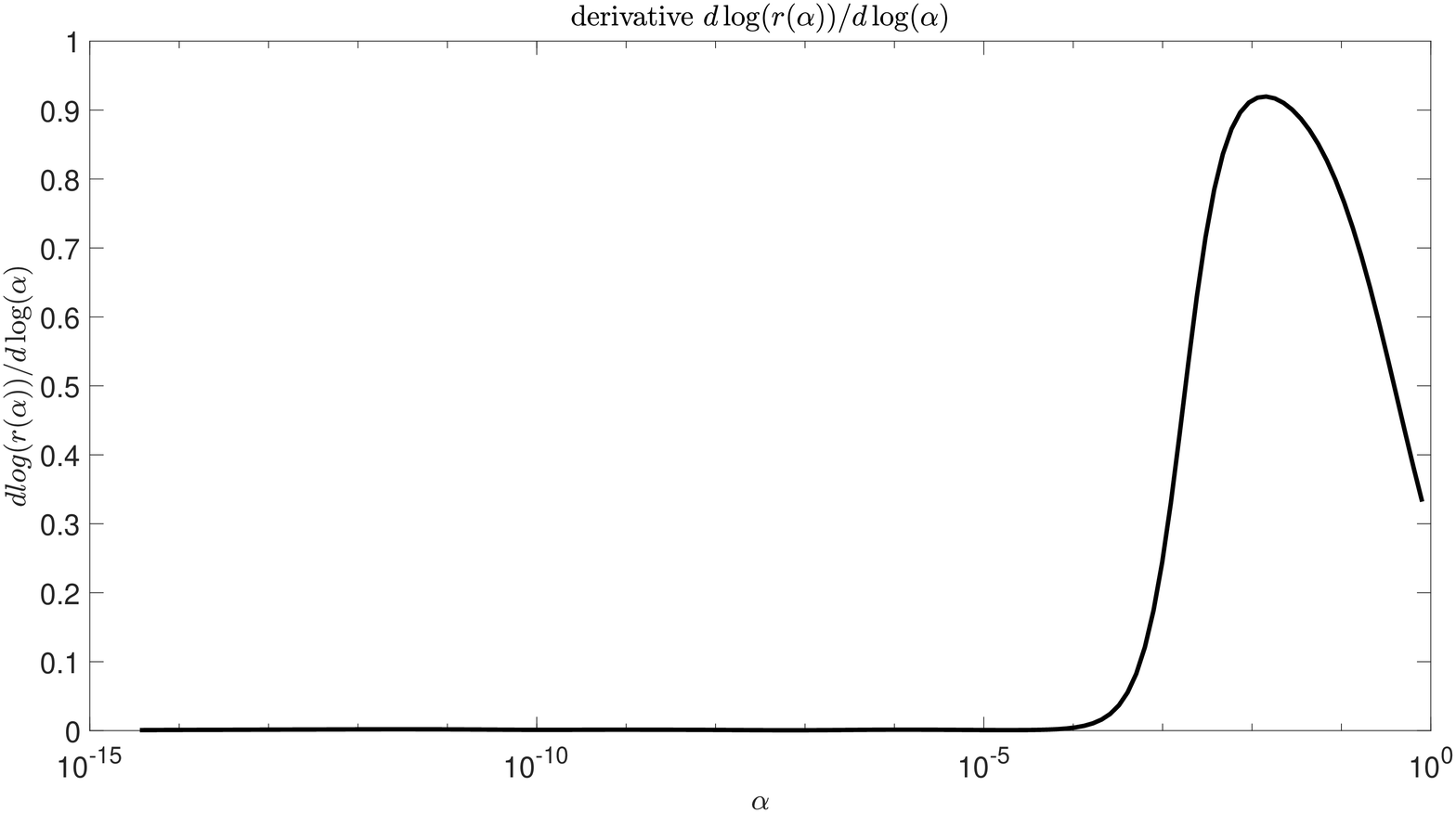}
\caption{Numerical Results for the logarithmic smoothness setting, \eqref{eq:sc_general_sum} with $\varphi(t)=(-\log(t))^{-\frac{3}{2}}$. Top: residuals, bottom: corresponding derivative. For sufficiently small $\alpha$ in the noise free case (left) we see the asymptotics $\|A\xa-A\xd\|\leq C \sqrt{\alpha}(-\log(\alpha))^{-\frac{3}{2}}$ shown in \cite{HofPlat}. Under noisy data this is not visible, only a predominantly concave part remains for the larger $\alpha$.}\label{fig:logsmooth}
\end{figure}

The opposite appears to be the case in the low smoothness setting, i.e. when $\varphi$ in \eqref{eq:sc_general} decays slower than a power function, where concavity appears to be the indicator. To demonstrate this, we consider a generalized source condition \eqref{eq:sc_general_sum} with $\varphi(t)=(-\log(t))^{-\kappa}$ for $\kappa=1.5$, see Figure \ref{fig:logsmooth}. It was shown in \cite{HofPlat} that $\|\xa-\xd\|\leq C (-\log(\delta)^{-\kappa})$. In this case one can further show that $\|A\xa-A\xd\|\leq C \sqrt{\alpha}(-\log(\alpha))^{-\kappa}$ for small enough $\alpha$. In the noise free case, this is plausible in the experiment. For noisy data this asymptotic is masked by the noise. What remains to extract information are the large regularization parameters. As usual, when $\alpha$ is too large, the residual changes little and then starts to drop. In the low smoothness case the residual curve is convex after the initial drop in the $\log$-$\log$-plot. To see this, we set 
\begin{align*}
R(\alpha):=&\log \|A\xa-y\|=\log\left(\sqrt{\alpha}(-\log(\alpha))^{-\kappa}\right)\\
=& \frac{1}{2}\log(\alpha)-\kappa\log(\log(\alpha)).
\end{align*}
Substituting $x:=\log(\alpha)$, we have $R(x)=\frac{1}{2}x-\kappa\log(x)$, and differentiating twice yields $\frac{d^2}{dx^2}R(x)=\frac{\kappa}{x^2}>0$ for all $x>0$, i.e., $r(x)$ is convex in the $\log$-$\log$-plot. This can also be seen in Figure \ref{fig:logsmooth}.

\begin{figure}
\includegraphics[width=0.49\linewidth]{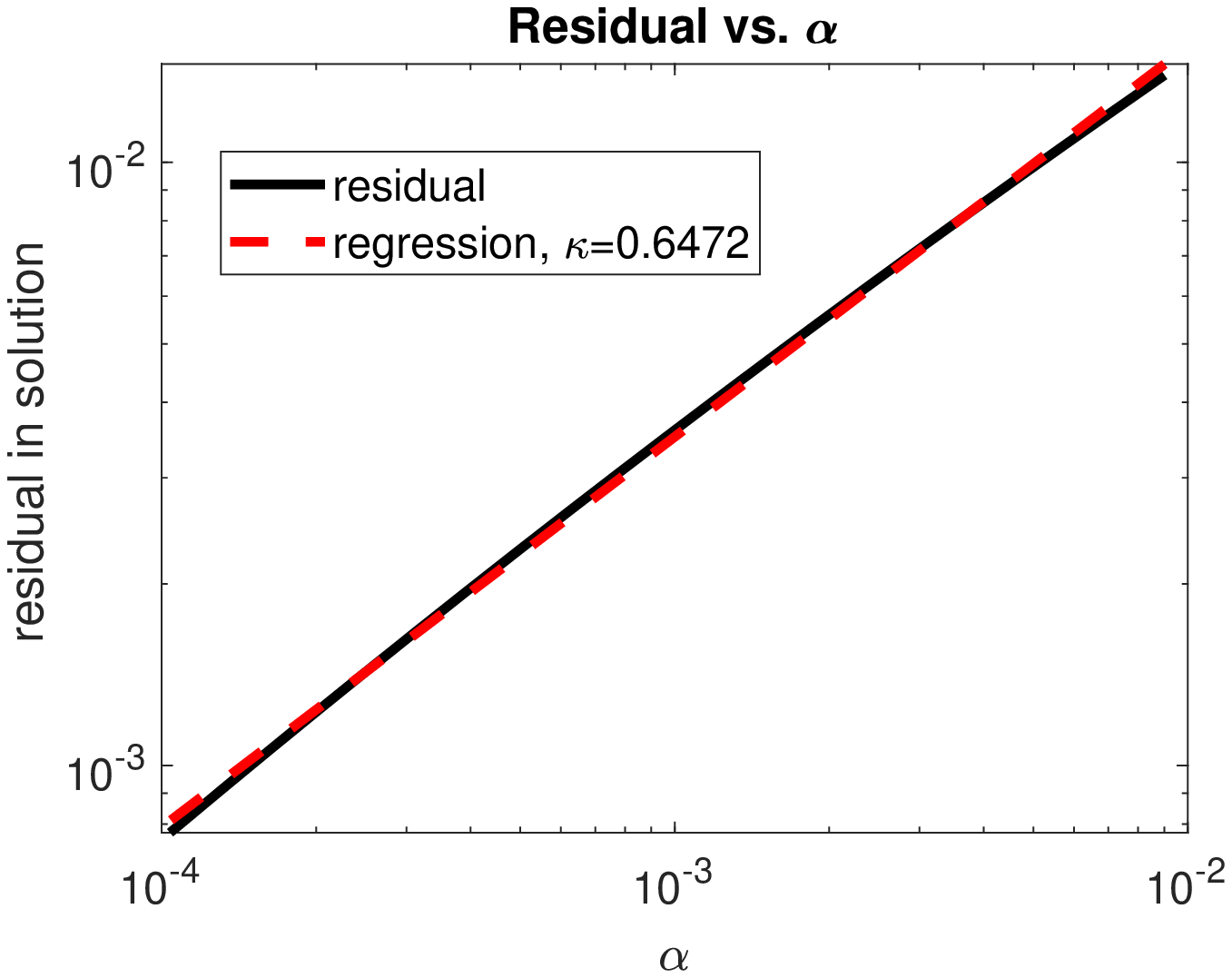}\includegraphics[width=0.49\linewidth]{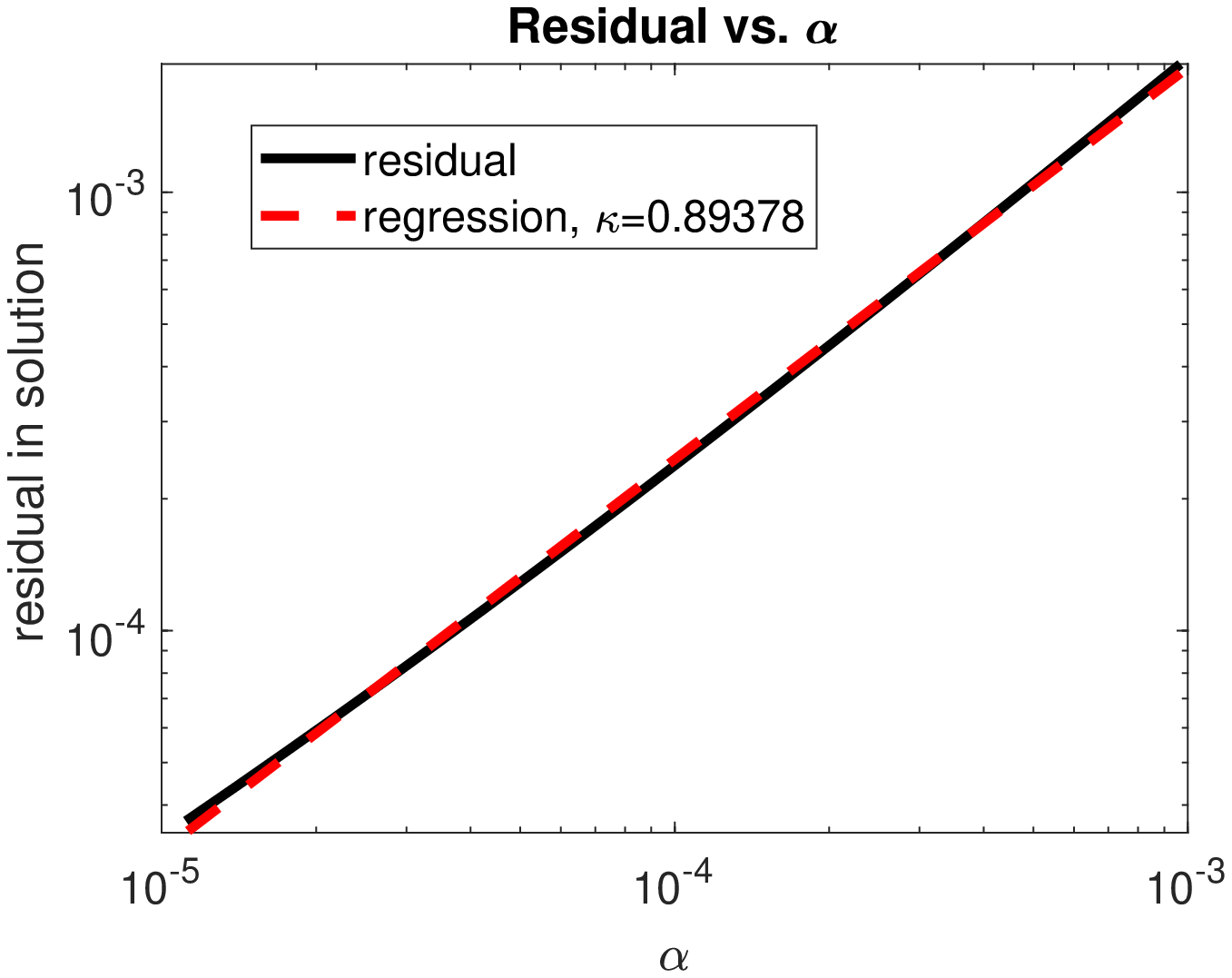}\caption{Zoom in on the residuals for exponential smoothness (left) and logarithmic smoothness (right). The deviation of the residual curves from the regression for the model $\|A\xa-y\|=c\alpha^\kappa$ is small. In the exponential case the residual is slightly concave, whereas in the logarithmic case it is slightly convex.}
\label{fig:lowhighzoom}
\end{figure}

The above observations indicate that it is possible to at least detect solution smoothness lower and higher than the H\"{o}lder type source condition \eqref{eq:sc_open} with $0<\mu<\frac{1}{2}$. In practice, however, this is difficult. Due to noise one has, in general, no access to the approximation rate $\|A\xa-y\|$ for sufficiently small $\alpha$ where the asymptotic rates become visible. Instead one will often be restricted to large regularization parameters, where a power-type regression often almost holds. Here one must find the minuscule differences and carefully inspect the deviation of the residual curve from the regression curve. This is exemplified in Figure \ref{fig:lowhighzoom}. With this observation we can also understand why high and low smoothness are difficult to handle in practice. In the range of regularization parameters one would expect in practice they behave almost like power-type source conditions.

\section{Experiments on Tomographic Data}\label{sec:tomo}
Moving away from the simulated data, we now apply the method to real data sets. This had been already done in \cite{DGLoja} using the Landweber algorithm to compute approximate solutions, but the results were difficult to interpret. We now apply our proposed Tikhonov-regularization approach. 

We use two samples from the tomographic X-ray data set collection provided by the Finnish Inverse Problems Society (FIPS), namely the data of a stuffed lotus root \cite{lotus} and the walnut data \cite{walnut}.

As for simulated data, we compute the approximate solutions by solving
\[
(A^T A+\alpha I)x=A^T y^\delta
\]
for various values of $\alpha$. The normal equation is solved using a Conjugate Gradient method. In particular we make use of the structure by interpreting $\alpha I$ as shifts of the matrix $A^T A$. This way only one Krylov subspace has to be built independent of the shifts $\alpha$, and for each fixed $\alpha$ the computations are cheap. This algorithm, which is described in detail in \cite[Algorithm 6]{FrommerMaass}, allows to use our method for high-dimensional problems and a large range of regularization parameters at low computational cost.

\begin{figure}
\includegraphics[width=0.99\linewidth]{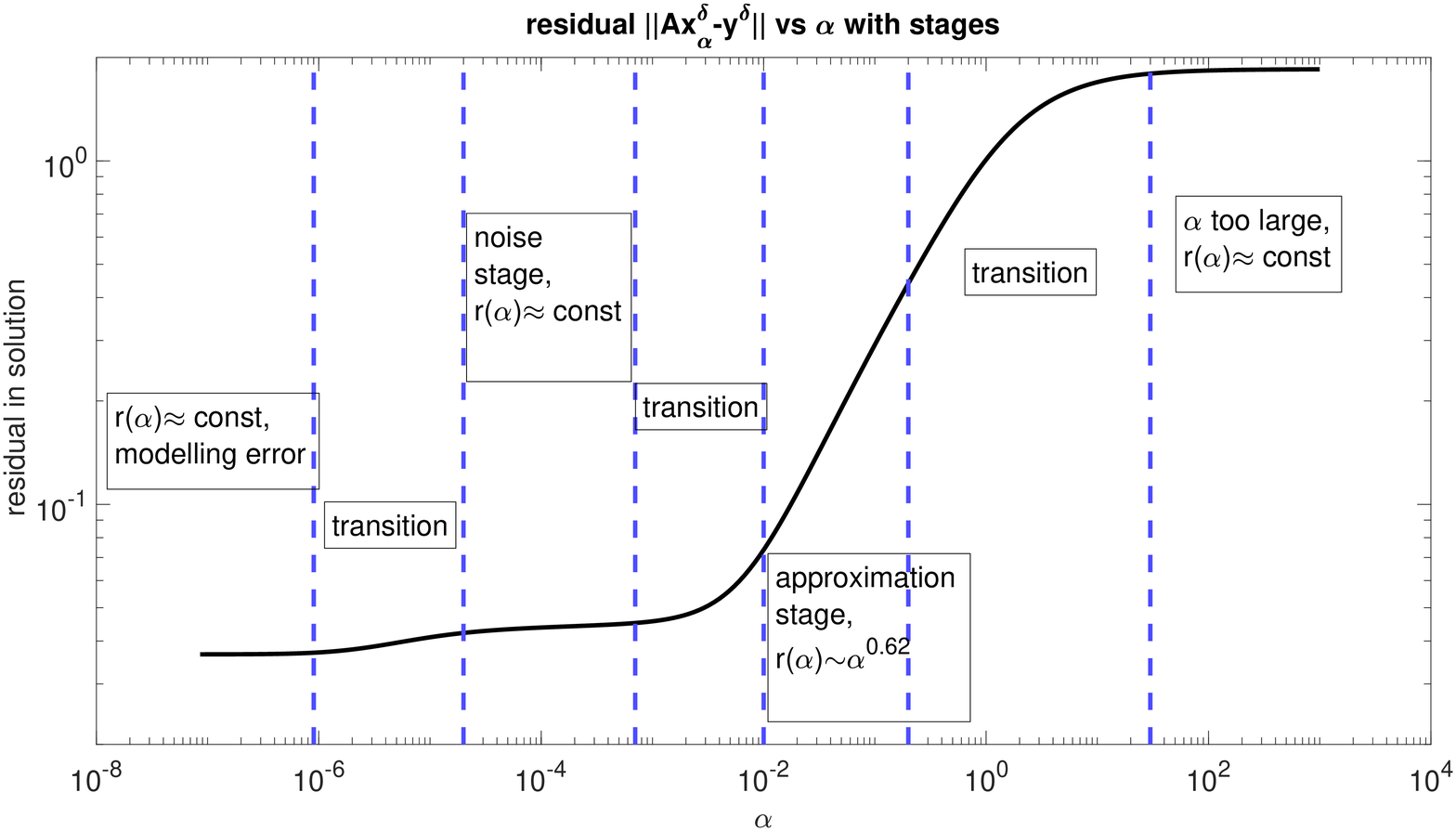}\\
\includegraphics[width=0.99\linewidth]{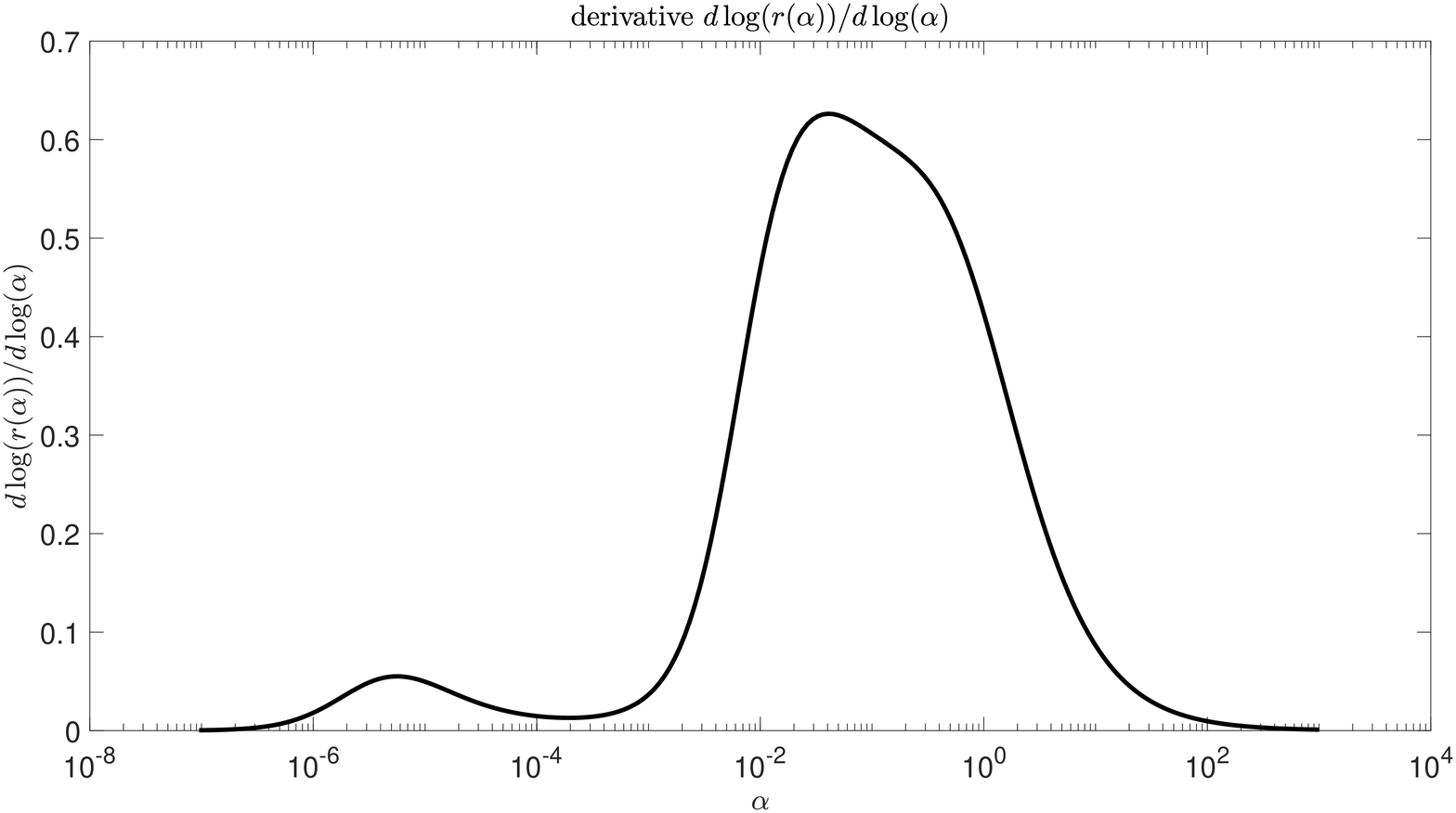}\caption{Residuals and its derivative for the large Lotus dataset. Except for the smallest $\alpha$, the plots have much similarity with the ones from simulated data, cf. Figure \ref{fig:residualcurve}, and we can observe the same phases. In particular, for $\alpha$ between approximately $0.01$ and $0.5$, we have the approximation phase that can be used to asses the solution smoothness. We zoom into this area in Figure \ref{fig:lotus_zoomedin}. For the smallest $\alpha$, we no longer have drop in residual but it stays constant, likely due to modelling errors. However, for $\alpha\sim 10^{-3}$ we observe the saddle point in the residual characteristic for the noise level. This residual with value $0.04$, yields the noise level $\delta\approx 0.04,$, or, together with $\|y^\delta\|=1.84$ and $\frac{\delta}{\|y^\delta\|-\delta}=\frac{0.04}{1.84-0.04}=0.022$, an estimated noise level of $2.2\%$.}\label{fig:lotus_zoomedout}
\end{figure}

\begin{figure}\includegraphics[width=0.99\linewidth]{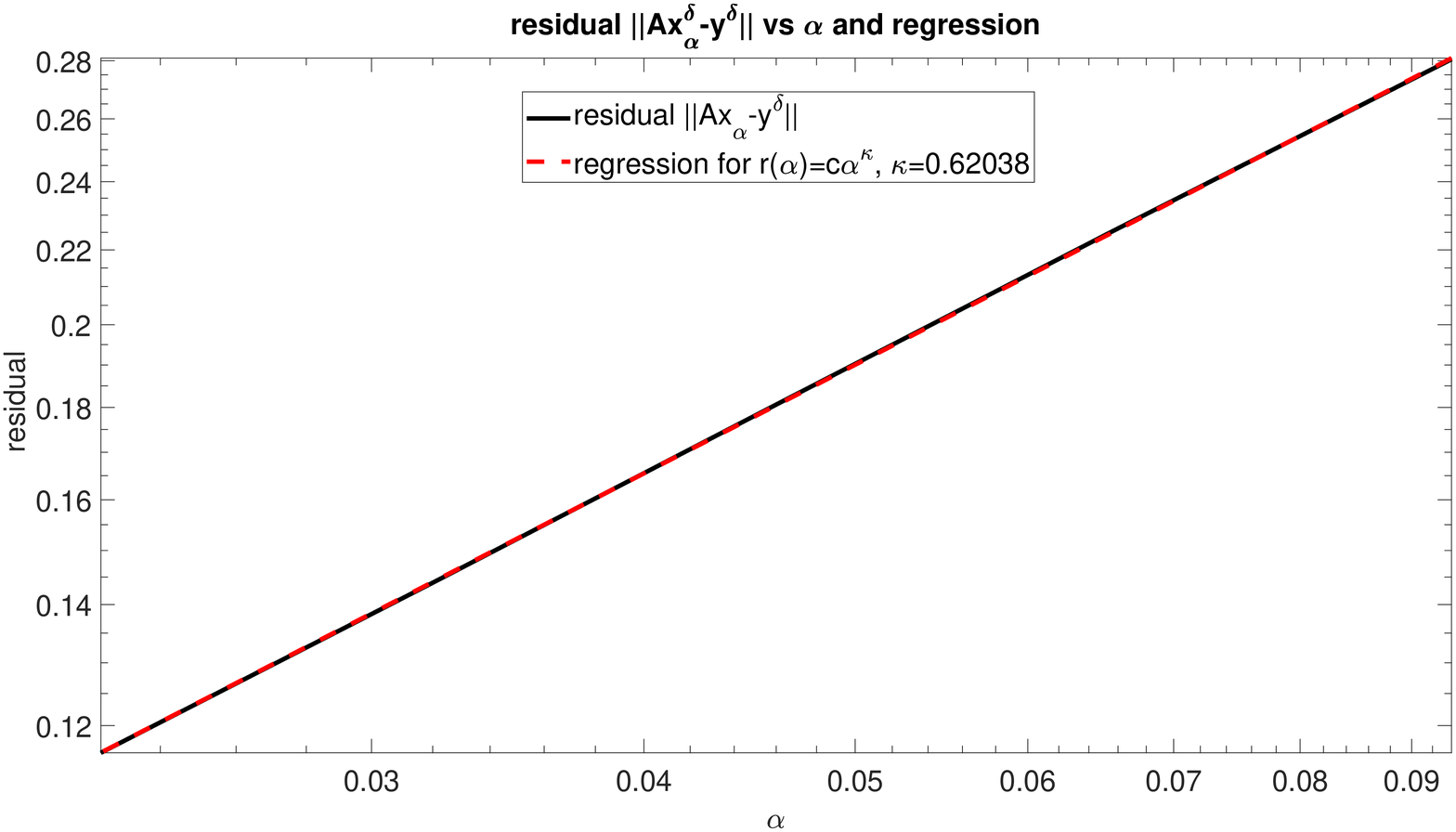}\caption{Results for the large Lotus dataset in the approximation phase, zoomed in from Figure \ref{fig:lotus_zoomedout}. Source element growth and residual are powers of $\alpha$, and the exponents correspond to a value $\mu\approx0.115$ in the source condition.}\label{fig:lotus_zoomedin}
\end{figure}

As our first example we consider the Lotus data set. For the largest data set \texttt{LotusData256.mat}, consisting of a matrix $A_{mn}\in \R^{51480\times65536}$ and measurements $y^\delta\in \R^{429\times 120}$. The results are shown in Figure \ref{fig:lotus_zoomedout}. We find that, as for simulated data, for large regularization parameters the residual is constant. In a second phase we find $\|A\xad-y^\delta\|\sim \alpha^{0.615}$ with the regression approach, see Figure \ref{fig:lotus_zoomedin} for a zoomed-in plot. This implies that $x^\dag$ fulfils a source-condition \eqref{eq:sc} with $\mu\approx 0.115$. Reducing $\alpha$ further yields a saddle-point like structure in the residual, from which we estimate the noise level $\delta\approx 0.04$ and a relative noise level of approximately $2.2\%$. After that, the residual drops slightly before it remains constant for all smaller regularization parameters, which we attribute to the modelling error.

\begin{figure}\includegraphics[width=0.99\linewidth]{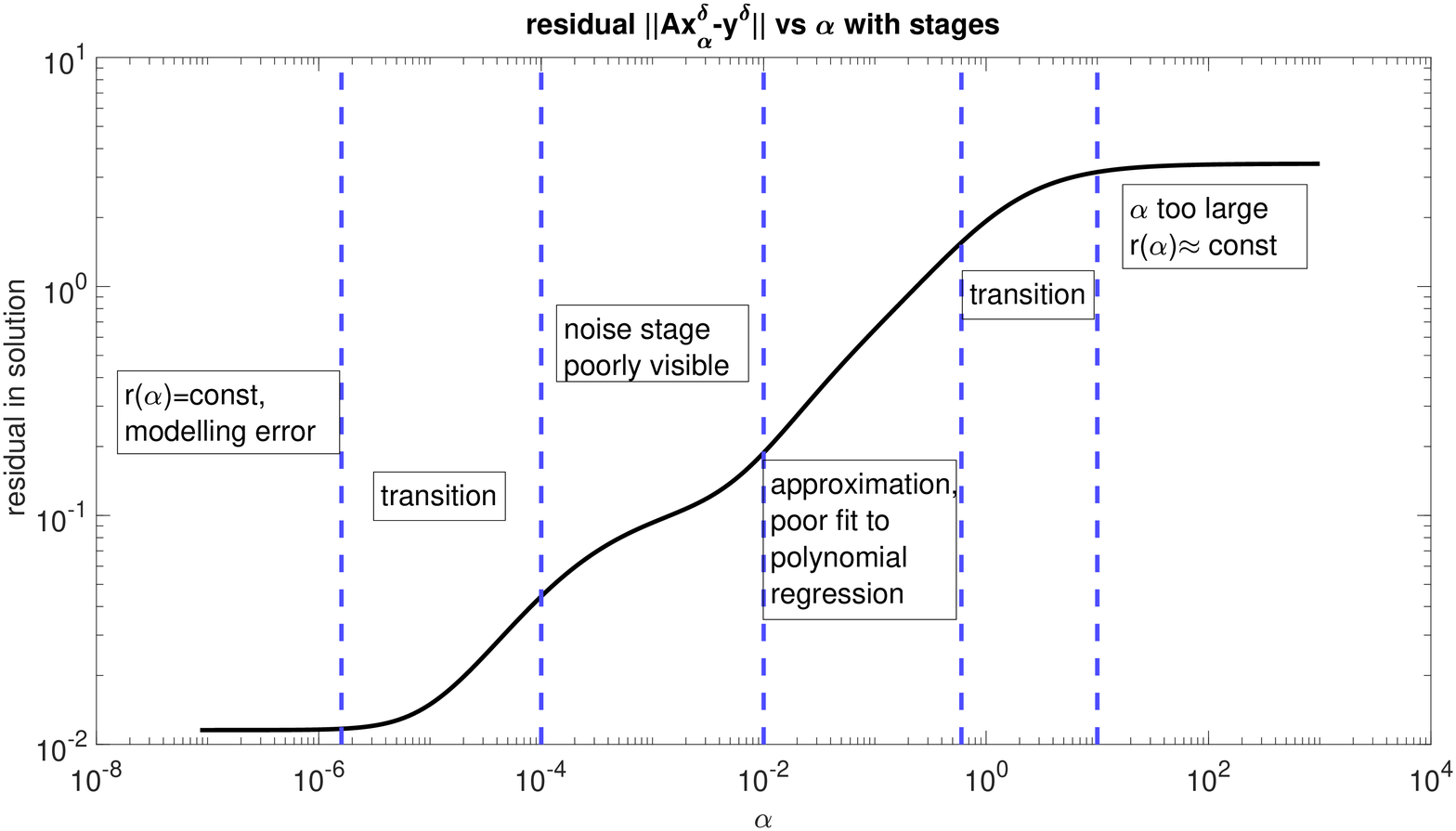}\\
\includegraphics[width=0.99\linewidth]{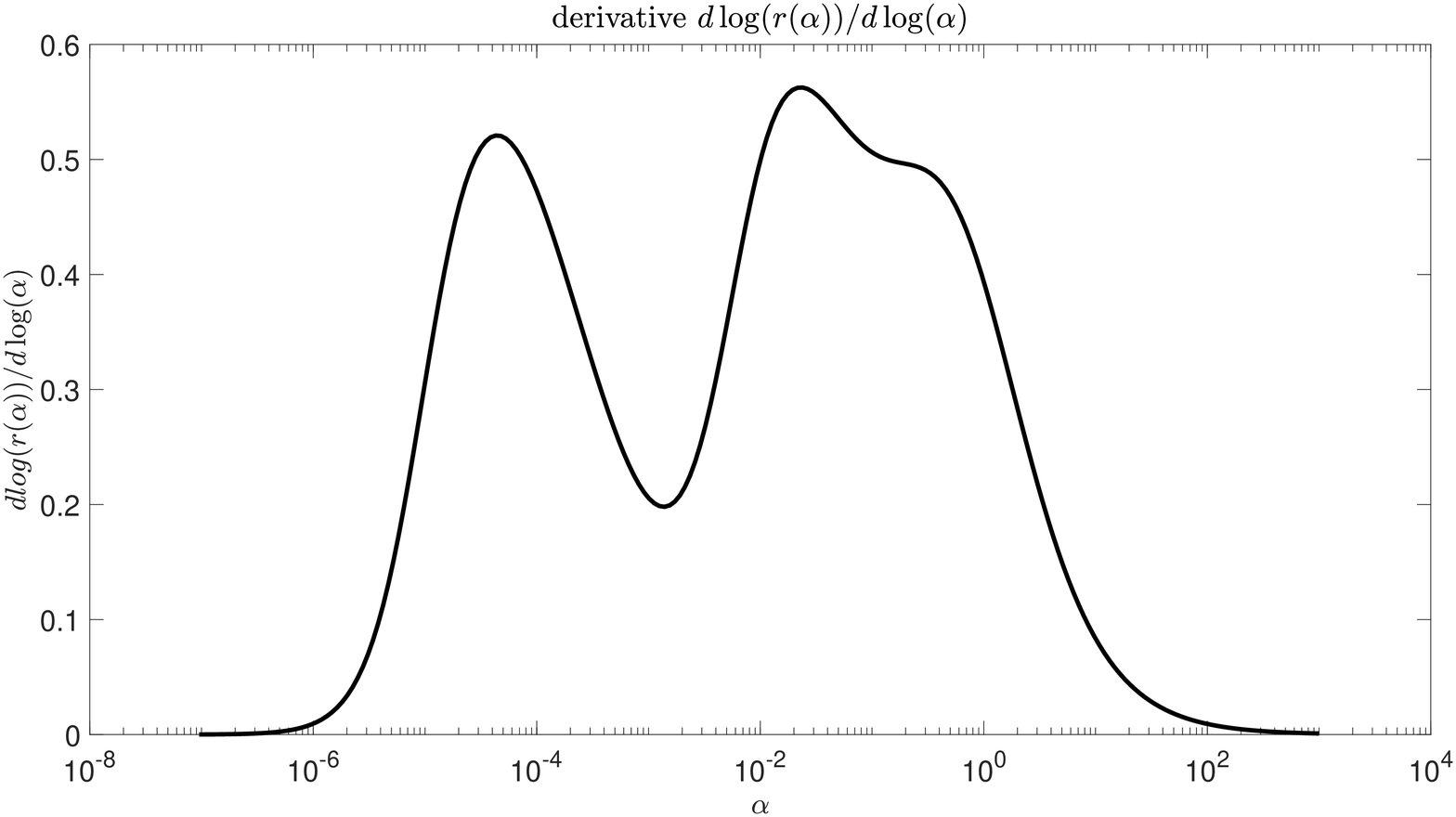}
\caption{Residuals and its derivative for the large walnut dataset.  After the burn in phase, we have a short phase the appears to be affine, but a closer inspecting in Figure \ref{fig:walnut_zoomedin} reveals that there are still oscillations around the linear regression. As with the Lotus data, we find that for the smallest $\alpha$ there is no more decrease in residual, likely for the same reason. }\label{fig:walnut_zoomedout}
\end{figure}

\begin{figure}\includegraphics[width=0.99\linewidth]{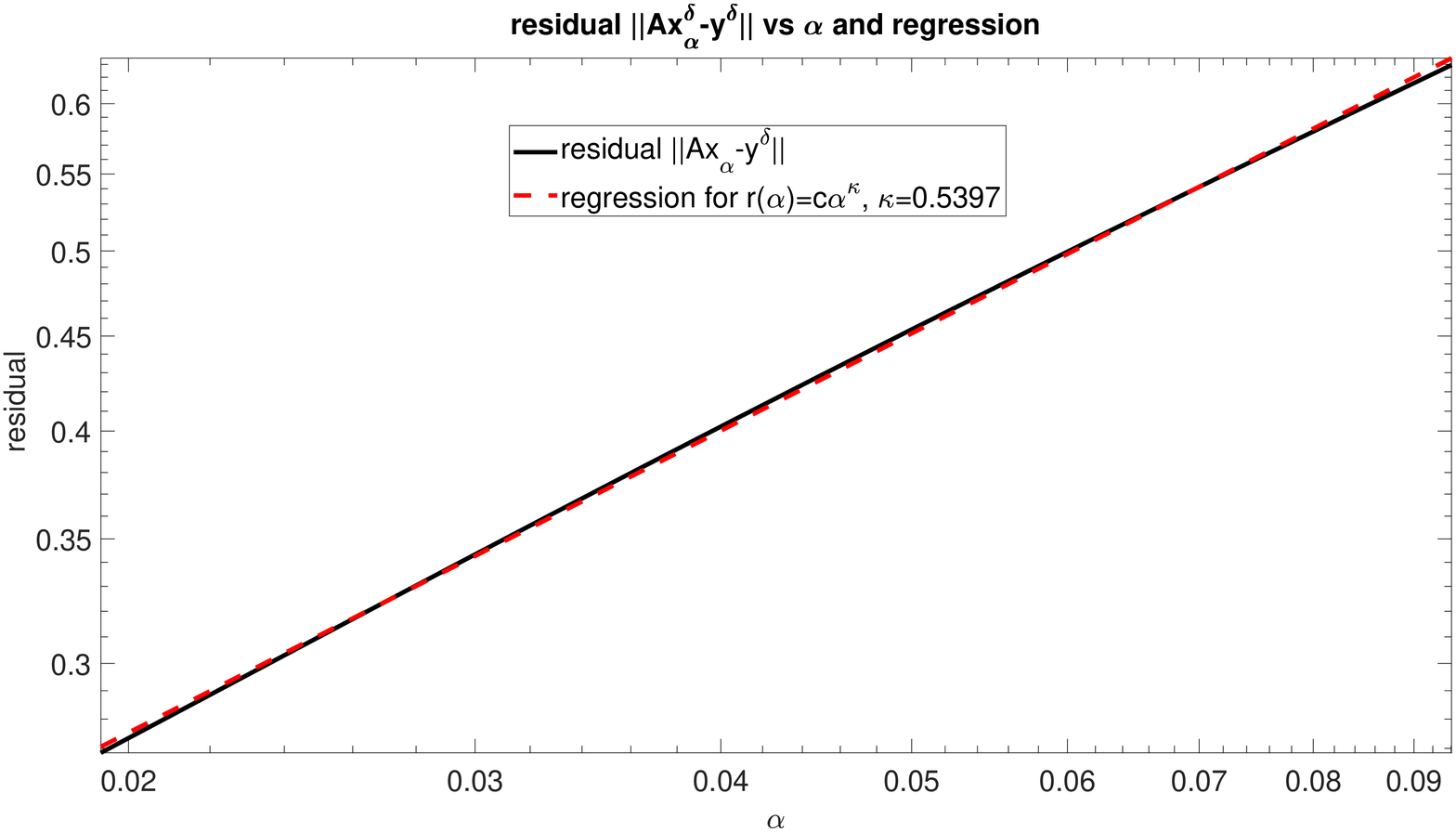}\caption{Residuals for the large walnut dataset in the approximation phase, zoomed in from Figure \ref{fig:walnut_zoomedout}. While the regression suggests a value $\mu\approx 0.02$, it can be seen that the residual oscillates around the regression lines, indicating the the source condition \eqref{eq:sc} does not hold for the walnut data.}\label{fig:walnut_zoomedin}
\end{figure}

We repeat the experiment for the largest data set of the walnut, see Figures \ref{fig:walnut_zoomedout} and \ref{fig:walnut_zoomedin}. In the approximation phase we find through regression that $\|A\xad-y^\delta\|\approx \alpha^{0.52}$, which would correspond to a source smoothness with $\mu=0.02$ . However, we see that the residual curve oscillates around the regression line, which would not be the case if a H\"{o}lder source condition would hold. Hence, we conclude that the walnut data does not fulfil a H\"{o}lder-type source condition \eqref{eq:sc_open} for any $\mu>0$.

\section{Choice of the regularization parameter}\label{sec:paramchoice}
A main task for the regularization of inverse problems in general and Tikhonov regularization in particular is the choice of the regularization parameter. Over time, so many parameter choice rules have been proposed such that it is not easy to keep a full overview. This is also not the purpose of this section. We will, on one hand, return to the optimization task \eqref{eq:minalpharescurve}, which we had introduced to estimate the noise level, and view it as a parameter choice rule. On the other hand, we will discuss the consequences of Theorem \eqref{thm:conv} for analysis of parameter choice rules in general. 

\subsection{Relating parameter choice rules}
Parameter choice rules are usually divided into three subgroups: a-priori choices $\alpha=\alpha(\delta)$, a-posteriori choices $\alpha=\alpha(\delta,y^\delta)$, and heuristics $\alpha=\alpha(y^\delta)$. A-priori choices, for Tikhonov regularization \eqref{eq:apriori}, are often regarded as the theoretical optimum that is infeasible in practice due to the lack of the necessary parameters $\mu$ and $\delta$. We have described above how both can be estimated, such that the a-priori choice can be carried out. However, as we have also noted, in particular the estimation of $\mu$ can be difficult. Therefore, we discuss in the following the impact of Theorem \ref{thm:conv} on other parameter choice rules.

The (Morozov) discrepancy principle is the most prominent a-posteriori principle. It only requires the knowledge of $\delta$ to select 
\[
\alpha^\ast=\sup\{ \alpha: \|A\xad-y^\delta\|\leq \tau \delta\}
\]
for some $\tau>1$. Of course, this requires the knowledge of $\delta$. Technically one also requires $\mu<\frac{1}{2}$ in the source conditions \eqref{eq:sc_sum} due to the restrictions of the residual as shown in Theorem \ref{thm:conv}. Heuristic parameter choice rules are often the only alternative in practical situations when neither $\mu$ nor the noise level are known. Naturally, most of these use the residual in some way or the other. Examples are the heuristic discrepancy principle , see. e.g., \cite{EHN}, where the functional
\begin{equation}\label{eq:heur_dp}
f(\alpha)=\frac{\|A\xad-y^\delta\|^2}{\alpha}
\end{equation} is minimized, or the L-curve \cite{hansen1993use}, and for the sake of brevity we will not discuss other methods.

To discuss the relation between the parameter choice rules, we start again with the a-priori choice \eqref{eq:apriori}, which is superior in the sense that it is applicable for $0< \mu\leq 1$ w.r.t a source condition \eqref{eq:sc} for $\xd$. Since the residual is non-informative for the interval $\frac{1}{2}<\mu<1$ we consider only $0<\mu<\frac{1}{2}$ in the following.
The a-priori choice is due to the decomposition
\begin{align*}
\|\xad-\xd\|^2&\leq \|\xa-\xd\|^2+\|\xad-\xa\|^2\\
&\leq c \alpha^{2\mu}+\frac{\delta^2}{\alpha},
\end{align*}
 see, e.g., \cite{EHN}, and optimizing over $\alpha$, i.e., minimizing
\begin{equation}\label{eq:apriori_heur}
f_{ap}(\alpha)=c \alpha^{2\mu}+\frac{\delta^2}{\alpha}.
\end{equation} 
Using the analogous decomposition for the residual, we have
\begin{align*}
\|A\xad-y^\delta\|^2&\leq \|A\xa-A\xd\|^2+\|A\xad-A\xa\|^2+\|A\xd-y^\delta\|^2\\
&\leq \tilde c \alpha^{2\mu+1}+\delta^2+\delta^2,
\end{align*}
where for the middle term the estimate
\begin{align*}
\|A\xad-A\xa\|^2&=\|A(\asa+\alpha I)^{-1} A^\ast(y-y^\delta)\|\\&\leq \|A(\asa+\alpha I)^{-1} A^\ast\|\,\|y-y^\delta\|\leq \delta
\end{align*}
is used. Inserting this into the heuristic discrepancy principle \eqref{eq:heur_dp} yields
\begin{align*}
f(\alpha)&=\frac{\|A\xad-y^\delta\|^2}{\alpha}\leq \tilde c\alpha^{2\mu} +2\frac{\delta^2}{\alpha},
\end{align*}
i.e., up to constants this is identical to the a priori functional \eqref{eq:apriori_heur}. Hence, asymptotically they yield the same convergence rate, and hence the heuristic discrepancy principle yields the optimal convergence rate $\|\xad-\xd\|\leq c\delta^{\frac{2\mu}{2\mu+1}}$ whenever $0<\mu<\frac{1}{2}$. Of course this convergence result is well-known (see, e.g., \cite{EHN}), but this proof is simpler than the standard argument.

For Morozovs discrepancy principle we use, for simplicity, the slightly different version to choose $\alpha$ such that $\|A\xad-y^\delta\|=\tau\delta$ for $\tau>1$. Considering that when $\|A\xad-y^\delta\|>\delta$ we have shown that $\|A\xad-y^\delta\|\approx \|A\xa-A\xd\|=\bigo(\alpha^{\mu+\frac{1}{2}})$, the discrepancy principle corresponds approximately to solving
\[
\tau \delta=\|A\xa-A\xd\|=\tilde c \alpha^{\mu+\frac{1}{2}}
\]
with solution $\alpha^\ast= c(\tau,\mu) \delta^{\frac{2}{2\mu+1}}$. This is, up to constants, the solution to the optimization problems for the a-priori choice \eqref{eq:apriori_heur} and the heuristic discrepancy principle \eqref{eq:heur_dp}. Therefore, all 3 methods differ only in constants and yield the same convergence rate. However, due to the differing constants, they do not yield the same regularization parameters in practice.

We now discuss the L-curve as one of the most popular heuristic parameter choice rules. For the L-curve method Tikhonov-approximations $\xad$ are calculated for many $\alpha$, similar to Algorithm \ref{alg}. The corresponding logarithms $\log(\|\xad\|)$ are then plotted against the logarithm of residuals $\log(\|A\xad-y^\delta\|)$. This typically yields a curve resembling the letter L, and the approximate solutions corresponding to the   corner of the L often are the closest to $\xd$. 

 As already noted e.g. in \cite{EHN}, the ``vertical'' part of the L is due to noise in the data, whereas the ``horizontal" part is due to approximation properties of the regularization method. More precisely, for Tikhonov regularization, we have that $\|\xad\|=\frac{\|A^\ast(A\xad-y^\delta)\|}{\alpha}$. As long as the residual is above the noise level, $\|A^\ast(A\xad-y^\delta)\|\approx \|A^\ast(A\xa-y)\|=\alpha$, which follows directly from the first-order condition for noise free data \eqref{eq:firstorder}. Therefore in this case $\|\xad\|\approx\mathrm{const}$. This yields the horizontal arm of the L. Once the noise level is reached, $\|A\xad-y^\delta\|\approx \delta$ for several magnitudes of $\alpha$ as explained in Section \ref{sec:noise}. Therefore, as $\alpha\rightarrow 0$, $\|\xad\|=\frac{\|A^\ast(A\xad-y^\delta)\|}{\alpha}\approx \frac{\|A^\ast(y-y^\delta)\|}{\alpha}\rightarrow \infty$ fast. This forms the vertical arm of the L.
The advantage of the L-curve over the other methods is that it essentially uses the gradient since $\|\xad\|=\frac{1}{\alpha}\|A^\ast(A\xad-y^\delta)\|$, and that the gradients have a distinct behaviour based on whether the residuals $A\xad-y^\delta$ are dominated by the approximation of $y$ or on the noise. Because for Tikhonov-regularization, the gradient does not depend on the smoothness of $\xd$, the L-curve does not require solution smoothness. Hence, it is more broadly applicable. On the other hand, the absence of smoothness information makes it more difficult to show convergence rates.

Nevertheless, if in particular $\xd$ satisfies a source condition, the corner of the $L$ corresponds to the regularization parameters where the residuals transitions from the asymptotic $\|A\xad-y^\delta\|\sim \alpha^{\mu+\frac{1}{2}}$ to $\|A\xad-y^\delta\|\approx\delta$. Therefore the L curve will yield a regularization parameter comparable to the discrepancy principle, and hence comparable to the other parameter choice rules as discussed above. Since the solution smoothness is irrelevant for the L-curve, more precise convergence results using our approach requires a more precise description of the noise, which is out of the scope of this paper.

\subsection{A new method and comparison}
The rule \eqref{eq:minalpharescurve} as crucial step in the estimation of the noise level can be used as a heuristic parameter choice rule. In our experiments, a slightly different variant appeared more favourable as the minimimum was less flat, and the functional resembled the curve for the reconstruction error much closer. Therefore, instead of \eqref{eq:minalpharescurve}, we search for the minimizer of
\begin{equation}\label{eq:resdiff}
f_{RD}(\alpha)=\argmin_\alpha \frac{\partial}{\partial \alpha}\|A\xad-y^\delta\|.
\end{equation}
We have shown in the previous sections that in the absence of noise the residual carries solution smoothness information and due to this behaves as $\|A\xa-y\|\sim\alpha^{\mu+\frac{1}{2}}$. As long as this is above the noise level, we also have $\|A\xad-y^\delta\|\sim \alpha^{\mu+\frac{1}{2}}$. We have also shown in Section \ref{sec:noise} that $\|A\xad-y^\delta\|\approx \delta$ for several magnitudes of $\alpha$ once it reaches the noise level. In between the two stages there is a transition phase. Therefore, the idea is to track the change in the residual when $\alpha$ decreases, which is done by \eqref{eq:resdiff}. We call this the residual differential method (RDM). A property separating it from other heuristic methods such as teh  heuristic discrepancy principle is that it tends to be optimistic, i.e., chooses regularization parameters lower than the optimal one, instead of being pessimistic, i.e., choosing too large regularization parameters. Therefore the approximate solutions appear slightly less smooth. Similarly to the L-curve, a convergence analysis for this method requires a more detailed description of the noise, which we will not pursue in this paper.

\begin{figure}
\includegraphics[width=\linewidth]{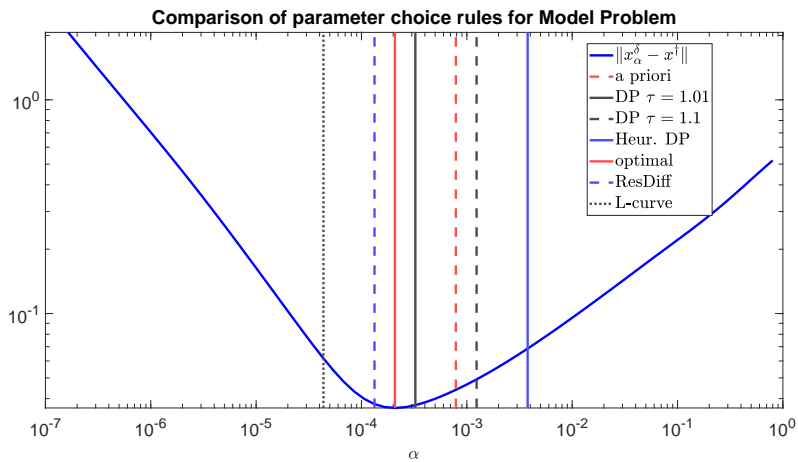}\caption{Comparison of different parameter choice rules for the Model Problem with $\eta=\beta=2$ and $0.5\%$ relative Gaussian noise.The values are representative but depend on the specific realization of the random noise.}\label{fig:pchoices}
\end{figure}

We compare some parameter choice rules in Figure \ref{fig:pchoices}. There we show the reconstruction errors $\|\xad-\xd\|$ for our Model Problem with $\eta=\beta=2$ ($\mu=0.375$ in \eqref{eq:sc_open}) for a large range of regularization parameters $\alpha$. $0.5\%$ Gaussian noise were added to the simulated exact data. We compare the following parameter choice rules: a-priori choice $\alpha=\left(\frac{\delta}{\rho}\right)^{\frac{2}{2\mu+1}}$ from \cite[Eq. (4.29)]{EHN}, where $\rho=\|w\|$ in the (numerical) source representation $\xd=\basa^\mu w$; Morozovs discrepancy principle with $\tau=1.01$ and $\tau=1.1$; the heuristic discrepancy principle; the practically infeasible parameter $\alpha_{opt}=\argmin_\alpha \|\xad-\xd\|$, the L-curve, and finally the RDM method \eqref{eq:resdiff}. The performance of the parameter choice rules is also compared in Table \ref{tab:pchoices}. Of course, the results are just a snapshot, and depending on the exact structure of the noise $\epsilon$ in \eqref{eq:noise} and the unknown solution $\xd$, the quality of the reconstructions for each parameter choice rule may vary slightly. The main point here was to demonstrate that, since they all work based on the change of the behaviour of the residual as function of $\alpha$ from predominantly approximating $y$ to being dominated by noise, they yield comparable results. In particular, this yields a rather simple explanation why heuristic parameter choice rules often do well in practice. A more detailed investigation and a relation to the theory of heuristic parameter choice rules (see, e.g., \cite{heur1,heur2,heur3}) is left as future work.
 
\begin{figure}
\includegraphics[width=\linewidth]{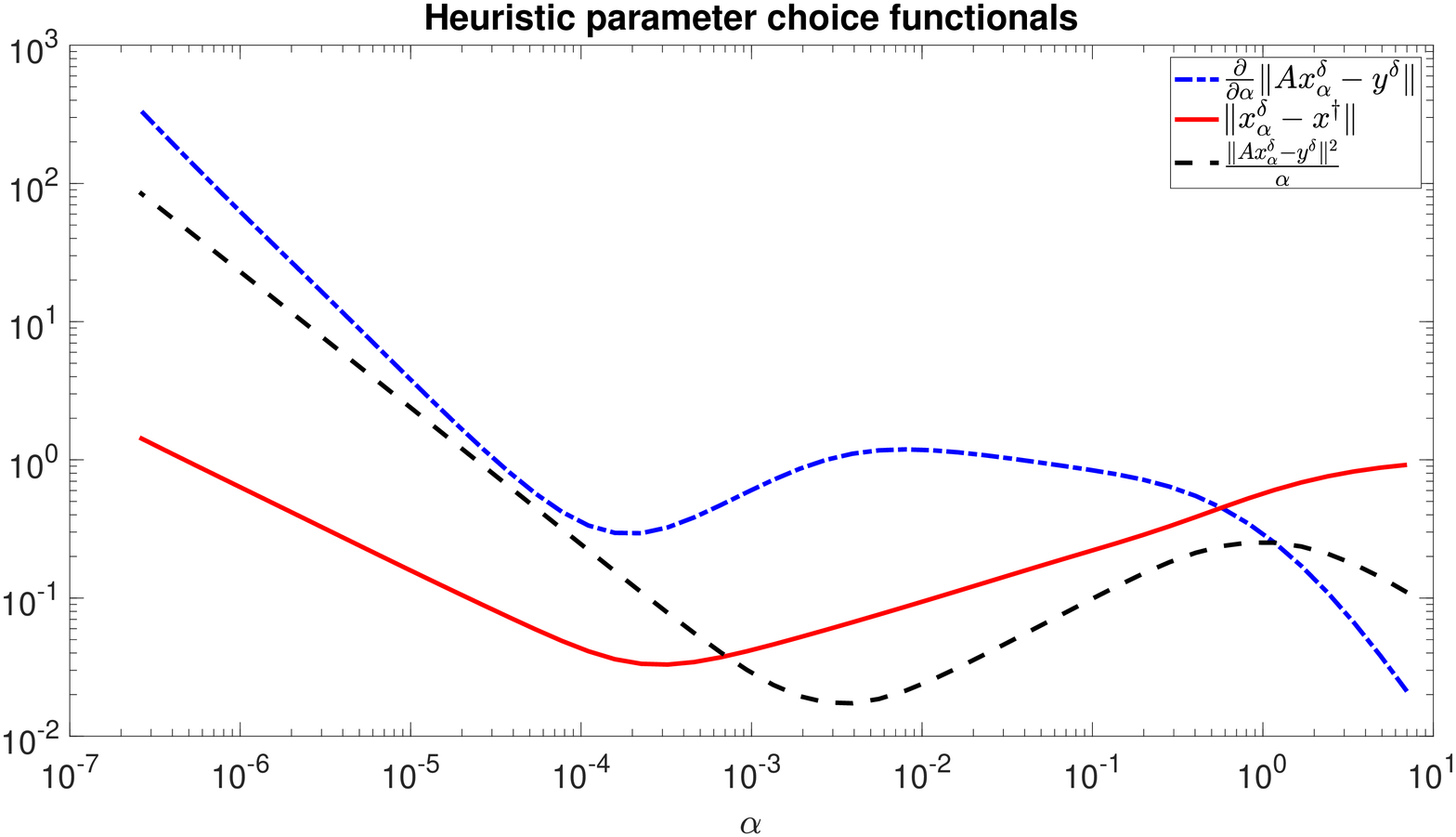}\caption{Comparison of the RDM functional with the reconstruction error and the heuristic discrepancy principle.}\label{fig:heuristics}
\end{figure}

\begin{table}\centering\caption{Comparison of different parameter choice rules for the Model Problem with $\eta=\beta=2$ and $0.5\%$ relative Gaussian noise. The values are representative but depend on the specific realization of the random noise.}
\begin{tabular}{c | c c c}
rule & $\alpha$ & $\frac{\|\xad-\xd\|}{\|x_{\alpha_{opt}}^\delta-\xd\|}$ & $\frac{\|A\xad-y^\delta\|}{\|Ax_{\alpha_{opt}}^\delta-y^\delta\|}$\\\hline
optimal &0.00032 & 1 & 1 \\
DP $\tau=1.01$ & 0.00046 & 1.04 & 1.1 \\
a priori & 0.00079 & 1.14 & 1.11 \\
DP $\tau=1.1$ & 0.0013 & 1.45 & 1.21 \\
Heur. DP & 0.0039 & 2.1 & 1.77\\
RDM & 0.00015 & 1.13 & 1.07
\end{tabular}\label{tab:pchoices}
\end{table}

\section{Conclusion and outlook}
We have shown that the classical source conditions $\xd\in \cR(\basa^\mu)$, $0<\mu<1$, are not only equivalent to approximation rates $\|\xa-\xd\|\sim \alpha^\mu$ under classical Tikhonov regularization, but also to a rate $\|A\xa-y\|\sim \alpha^{\frac{1}{2}+\mu}$ for $0<\mu<\frac{1}{2}$. This result allows to extract the solution smoothness from the residuals, making this information accessible in practical computation. We have demonstrated how larger and higher smoothness are detectable from the residual, although not quantifiable. We have demonstrated that an estimate of the noise level can be read off the residual curve. Because the residual carries so much information, we were able to relate parameter choice rules in a novel way. There are several open topics for future work. One possible topic is the extension of the smoothness estimation to modern Banach space regularization methods, which at first requires to find an appropriate way to measure solution smoothness, since source conditions are in general not applicable. A second path is to investigate the noise even further and to categorize different levels of noise smoothness. This could be combined with revisiting the heuristic parameter choice and clarify the open questions.

\section*{Acknowledgements}
D. Gerth was supported by Deutsche Forschungsgemeinschaft (DFG), project GE3171/1-1.  R. Ramlau was supported by the Austrian Science Fund (FWF, SFB ``Tomography across the scales'' F6805-N36.  D. Gerth would like to thank Prof. Oliver Ernst (TU Chemnitz) for the helpful comments.

\section*{Appendix}
We show 
\[
\sum_{i=1}^k \sigma_i^{q-2p}\langle\xd,v_i\rangle^2=\bigo(\sigma_k^{q+4\mu-2p})
\]
if $\sum_{i=k}^\infty \langle x,v_i\rangle^2=\bigo(\sigma_k^{4\mu})$ and $q+4\mu-2p<0$. 

For $k=1$ it is 
\begin{align*}
\sum_{i=1}^k \sigma_i^{q-2p}\langle\xd,v_i\rangle^2=\sigma_1^{q-2p}\left(\sum_{i=1}^\infty \langle x,v_i\rangle^2-\sum_{i=2}^\infty \langle x,v_i\rangle^2 \right)=\bigo(\sigma_i^{q+4\mu-2p}.
\end{align*}
Now assume the assertion holds for arbitrary $k>0$. Then
\begin{align*}
\sum_{i=1}^{k+1} \sigma_i^{q-2p}\langle\xd,v_i\rangle^2&=\sum_{i=1}^k \sigma_i^{q-2p}\langle\xd,v_i\rangle^2+\sigma_{k+1}^{q-2p}\langle x,v_i\rangle^2\\
&=\bigo(\sigma_k^{q+4\mu-2p})+\sigma_{k+1}^{q-2p}\left(\sum_{i=k+1}^\infty \langle x,v_i\rangle^2-\sum_{i=k+2}^\infty \langle x,v_i\rangle^2 \right)\\
&=\bigo(\sigma_k^{q+4\mu-2p})+\bigo(\sigma_{k+1}^{q+4\mu-2p})=\bigo(\sigma_{k+1}^{q+4\mu-2p})
\end{align*}
because $q+4\mu-2p<0$.

\bibliographystyle{plain}

\end{document}